\UseRawInputEncoding

\documentclass[11pt]{amsart}
\usepackage[utf8]{inputenc}
\usepackage[toc,page]{appendix}

\usepackage[numbers]{natbib}
 \usepackage{tikz-cd}
 \usepackage{comment}
\usepackage[margin=1in]{geometry}
\usepackage{amsmath,amsthm,amssymb}

\usepackage{amsthm}

\usepackage[hidelinks]{hyperref}
\usepackage[capitalise,nosort,nameinlink]{cleveref}

\newtheorem{theorem}{Theorem}[section]

\newtheorem{corollary}[theorem]{Corollary}

\newtheorem{proposition}[theorem]{Proposition}
\newtheorem{lemma}[theorem]{Lemma}

\theoremstyle{definition}
\newtheorem{definition}[theorem]{Definition}
\newtheorem{example}[theorem]{Example}
\newtheorem{notation}[theorem]{Notation}
\newtheorem{construction}[theorem]{Construction}

\theoremstyle{remark}
\newtheorem{remark}[theorem]{Remark}

\crefname{theorem}{Theorem}{Theorems}
\crefname{corollary}{Corollary}{Corollaries}
\crefname{proposition}{Proposition}{Propositions}
\crefname{lemma}{Lemma}{Lemmas}

\crefname{definition}{Definition}{Definitions}
\crefname{example}{Example}{Examples}
\crefname{notation}{Notation}{Notations}
\crefname{description}{Description}{Descriptions}
\crefname{construction}{Construction}{Constructions}

\crefname{remark}{Remark}{Remarks}

\usepackage{floatrow}
\usepackage[textwidth=2cm,textsize=small,colorinlistoftodos]{todonotes} 

\newcommand{\sSet}{\mathsf{sSet}}
\newcommand{\Pos}{\mathsf{Pos}}
\newcommand{\Sd}{\mathrm{Sd}}
\newcommand{\Ex}{\mathrm{Ex}}
\newcommand{\Set}{\mathsf{Set}}
\newcommand{\Cat}{\mathsf{Cat}}
\newcommand{\sSetBC}{\mathsf{sSet_{un}}}
\newcommand{\sSetB}{\mathsf{sSet_{ns}}}

\newcommand{\Osc}{\mathsf{Osc}}

\newcommand{\OrdSet}{\mathsf{OrdSet}}

\newcommand{\Hom}{\mathrm{Hom}}
\newcommand{\colim}{\mathrm{colim}}

\newcommand{\id}{\mathrm{id}}

\begin{document}
 
 
\title{Homotopy theory for ordered simplicial complexes}
\author{Melissa Wei} 
\address{School of Mathematics, University of Minnesota, Minneapolis MN, 55415, USA}
\email{wei00359@umn.edu}
\thanks{This material is based upon work supported by the National Science Foundation under Award No. DMS-2528362. }
 \begin{abstract}

We construct a model structure on the category of ordered simplicial complexes, Quillen equivalent to the standard model structure on simplicial sets. This shows that simplicial complexes, which are fully combinatorial in nature, provide a new model for the homotopy theory of spaces.
\end{abstract}
\maketitle

\setcounter{tocdepth}{1}
\tableofcontents
\section{Introduction}
Abstract simplicial complexes are completely combinatorial objects which are simple to define, and their geometric realizations are the familiar geometric simplicial complexes \cite{FriedmanSimplicialSets}. As such, simplicial complexes are easy to understand and to work with. Beyond topology, where they give rise to the notion of the triangulation of a space, simplicial complexes are used in many applications in other areas of mathematics, such as sphere packing combinatorics \cite{HalesSpherePacking1992}, the Whitney complex of a graph \cite{LarrionPizanaVillarroelClique},  and Vietoris--Rips complexes \cite{ChambersDeSilvaEricksonRips}; we can even form the Stanley--Reisner ring from simplicial complexes, giving a link to commutative algebra \cite{FranciscoMerminSchweigSR}. However, there is a question of how well they are able to approximate topological spaces. In particular, are all spaces homeomorphic to simplicial complexes, or in other words, can all spaces be triangulated? The disproof of the triangulation conjecture \cite{manolescu2015pin2equivariantseibergwittenfloerhomology} shows this is false for even very well-behaved spaces such as manifolds. 

On the other hand, simplicial sets can model all topological spaces up to homotopy. Simplicial sets are more flexible than simplicial complexes which has made them ideal for this purpose. Therefore, a natural question to ask would be whether simplicial complexes can also model topological spaces up to homotopy. We answer this in the affirmative in the case of ordered simplicial complexes.

To do this, we establish the homotopy theory of two categories whose objects represent simplicial complexes: ordered simplicial complexes (denoted $\Osc$), and simplicial sets characterized by the property that a set of vertices has at most one nondegenerate simplex with exactly those vertices (denoted $\sSetBC$); the latter is precisely the subcategory of simplicial sets with Properties B and C as described in Definition 12.1.1 of \cite{finitespaces}.  Previously,  as in Theorem 12.1.8 of \cite{finitespaces}, ordered simplicial complexes were thought to embed fully faithfully into simplicial sets, and were thus equivalent to $\sSetBC$. However, this is untrue, as the existence of ordered loops which are valid simplicial sets in $\sSetBC$, but do not constitute ordered simplicial complexes, make the relation between these two categories significantly more subtle (see Remark \ref{rmkloop}). 

Note that we consider ordered simplicial complexes as opposed to unordered simplicial complexes. The category $\Osc$ is more well-behaved as opposed to the category of simplicial complexes; for instance, products are preserved by geometric realization in $\Osc$ but not in simplicial complexes (see \cite{bergner2024simplicialsetstopologycategory}, Example 3.2). Note that there is no loss of topological generality that arises from this restriction; in particular, there is a canonical barycentric subdivision functor from simplicial complexes to ordered simplicial complexes which gives the same geometric realization (described in Chapter 10.1 of \cite{finitespaces}). Also, a predefined order is necessary for a natural passage of simplicial complexes to simplicial sets.

To establish a model structure on $\Osc$, we need two intermediate categories between $\Osc$ and $\sSet$: $\sSetB$, the category of nonsingular simplicial sets, which was studied in \cite{fjellbonsSet}, and $\sSetBC$. 
Ordered simplicial complexes also have the property that a set of $n$ vertices can only have one $n$-simplex with those vertices, but the ordering on vertices must be compatible between simplices. Thus, we take the transitive closure of the underlying graphs of simplicial sets in $\sSetBC$ in order to lift its model structure to $\Osc$. None of these are equivalent as categories, but they are Quillen equivalent. Hence, both categories, which encode simplicial complexes in different ways, are suitable models for the homotopy theory of spaces. We prove the following:
\begin{theorem}
The category $\sSetBC$ admits a model structure Quillen equivalent to the standard model structure on $\sSet$ and to the model structure on $\sSetB$ of \cite{fjellbonsSet}.
\end{theorem}
For the full statement, see Theorem \ref{sSetBCmodel}. Specifically, the model structure of $\sSetB$, is lifted along the Quillen adjunction \[\begin{tikzcd}
            \sSetB \arrow[r, shift left=.75ex,"L"{name=G}] & \sSetBC\arrow[l, shift left=.75ex, "i"{name=F}] 
            \arrow[phantom, from=F, to=G, "\scriptstyle\dashv" rotate=270]     
        \end{tikzcd} \]to a cofibrantly generated model structure.

This is an intermediate step to the next main result of our paper, in which we prove that $\Osc$ has a model structure which is Quillen equivalent to $\sSet$:

\begin{theorem}
There is a cofibrantly generated model structure on $\Osc$ transferred from $\sSetBC$ via the Quillen pair \[\begin{tikzcd}
            \sSetBC \arrow[r, shift left=.75ex,"F"{name=G}] & \Osc\arrow[l, shift left=.75ex, "U"{name=F}] 
            \arrow[phantom, from=F, to=G, "\scriptstyle\dashv" rotate=270]      
        \end{tikzcd} \]
        Furthermore, this adjunction is a Quillen equivalence, and hence $\Osc$ is Quillen equivalent to $\sSet$. 
\end{theorem}
See Theorem \ref{maintheorem2} for more details.
An illustration of this chain of equivalences is below:
\[
        \begin{tikzcd}
            \sSet \arrow[r, shift left=.75ex,"D\Sd^2"{name=G}] & \sSetB\arrow[l, shift left=.75ex, "\Ex^2i"{name=F}] 
            \arrow[phantom, from=F, to=G, "\scriptstyle\dashv" rotate=270]\arrow[r, shift left=.75ex,"L"{name=H}] & \sSetBC\arrow[l, shift left=.75ex, "i"{name=I}] 
            \arrow[phantom, from=H, to=I, "\scriptstyle\dashv" rotate=270] \arrow[r, shift left=.75ex,"F"{name=J}] & \Osc\arrow[l, shift left=.75ex, "U"{name=K}] 
            \arrow[phantom, from=J, to=K, "\scriptstyle\dashv" rotate=270]   
        \end{tikzcd}
\]
where the first Quillen equivalence was established in \cite{fjellbonsSet}. We establish Quillen equivalences along the right two adjunctions to complete the passage from $\sSet$ to $\Osc$.

As a note to the reader, there have been independently developed efforts to establish a model structure on simplicial complexes in \cite{minichiello2025thomasontypemodelstructuressimplicial}, which is Quillen equivalent to simplicial sets. This model category of unordered simplicial complexes is referred to as $\mathsf{Cpx}$, but it is fundamentally different from both model categories $\sSetBC$ and $\Osc$ we establish here. In particular, $\mathsf{Cpx}$ is not a subcategory of simplicial sets, and though \cite{minichiello2025thomasontypemodelstructuressimplicial} constructs a fully faithful functor $Sing\colon  \mathsf{Cpx} \to \sSet$, it does not even land in $\sSetBC$. Thus its essential image does not give simplicial sets with the expected properties of simplicial complexes.

\subsection*{Acknowledgments}
Thank you so much to Inna Zakharevich and Brandon Shapiro, and also to my advisors Maru Sarazola and Mike Hill, for all their help and guidance. This project would not have been possible without their support, suggestions, and mentorship. I would also like to thank Samuele Pedretti for bringing an error in an earlier version of Section 4 to my attention, and to the anonymous referee who suggested looking into properties of $\sSetBC$ and $\Osc$, which resulted in the last section of this paper. 
\section{Preliminaries and Properties of Simplicial Sets}
In this section, we give background on  relevant properties of simplicial sets, especially those which make them resemble simplicial complexes. We will also review a combinatorial framework for working with subdivision in simplicial sets. This will be used to relate certain simplicial sets to posets, which will be relevant in Section 5.

We start by reviewing the definition of simplicial complexes:
\begin{definition}
    An \textit{(unordered) simplicial complex} $(S$, $K \subseteq \mathcal{P}(S))$ is a set of vertices $S$ and a collection of subsets $K$ such that: (1) $\forall v \in S, \{v\} \in K$, and (2) if $T \in K$ and $U \subseteq T$ then $U \in K$.
   
\end{definition}
In particular, given a set of vertices in a simplicial complex, there is at most one simplex with exactly those vertices. 

We return to the category $\sSet$. A subcategory of simplicial sets which captures the condition that nondegenerate simplices have distinct vertices is that of nonsingular simplicial sets (defined from Definition 1.1 of \cite{fjellbonsSet}), which are a full subcategory of simplicial sets whose objects are defined as follows. 

\begin{definition}
A simplicial set $X$ is \textit{nonsingular} if for every nondegenerate simplex $x$, the representative map $\bar{x}$ of the correspondence $X_n \cong \sSet(\Delta[n], X), x \mapsto \bar{x}$ under the Yoneda lemma is injective. We will sometimes refer to such a simplex $x$ as $\textit{embedded}$. The category of nonsingular simplicial sets will be denoted $\sSetB$.
\end{definition}
Nonsingular simplicial sets, where each isolated simplex cannot have repeated vertices, are essentially an intermediate class of simplicial sets in the passage towards simplicial complexes, where each such simplex is also unique.

Note that nonsingular simplicial sets are not necessarily simplicial complexes. Although both have the property that none of the simplices have repeated vertices, nonsingular simplicial sets may have multiple simplices with exactly the same distinct vertices, while simplicial complexes can have at most one such simplex.

Now we reference several properties of simplicial sets:  Properties B and C are defined in \cite[page 85]{finitespaces}.
\begin{definition} Given a simplicial set $X$, we define the following properties: 
\label{props}
\begin{itemize}
\itemsep-1em 
\item \textit{Property B} (the distinct vertex property): Every nondegenerate simplex of $X$ has distinct vertices $\\$
\item \textit{Property C} (the unique simplex property): Given a set of distinct vertices, there is at most one nondegenerate simplex with those vertices. $\\$
\end{itemize}

\end{definition}
By the end of this section, we will have introduced three categories of simplicial sets: $\sSet$, $\sSetB$, and $\sSetBC$, which are related by a chain of strict inclusions. In particular, we will show $\sSetB$ is the category of simplicial sets with Property B and by definition $\sSetBC$ is the category of simplicial sets with both Properties B and C. In the figure below, we show some examples of simplicial sets in each category.

        \includegraphics[width=0.9\linewidth]{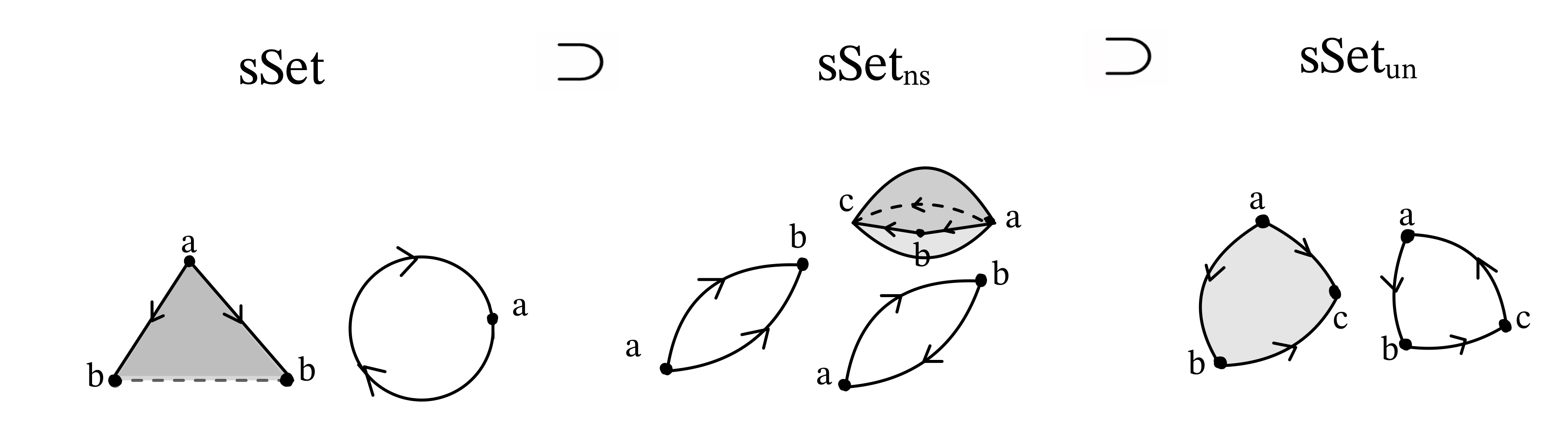}

In the rest of this section, we show how Properties B and C characterize some of the structures we wish to study.
\subsection{Property B}
\begin{proposition}
A simplicial set is nonsingular if and only if it has Property B.
\label{prop3}
\end{proposition}
\begin{proof}
If a simplicial set $X$ is nonsingular, then all of its simplices are embedded, so the representing map of each $x \in X_k$, $\bar{x}: \Delta[k] \rightarrow X$, is injective. When precomposing with the $k$ maps corresponding to each inclusion $[0] \rightarrow [k]$, we obtain $k$ distinct maps $\Delta[0] \rightarrow X$ representing the $k$ distinct vertices of $X$. Conversely, if $X$ is singular, then it must have a nonembedded, nondegenerate $k$-simplex $x$. Its representing map $\bar{x}: \Delta[k] \rightarrow X$ is not injective, so there must be two identical maps $\Delta[0] \rightarrow X$  obtained from precomposition of distinct maps $[0] \rightarrow [k]$ that are both vertices of $x$. Hence, $x$ cannot have $k$ distinct vertices, so $X$ cannot have Property B.
\end{proof}

Given any simplicial set, there is a canonical way to obtain a nonsingular simplicial set by applying \textit{desingularization}, as introduced in \cite[Remark 2.2]{waldhausen2013spaces}.

\begin{definition}

The \textit{desingularization} functor $D: \sSet \to \sSetB$ maps a simplicial set $X$ to the image of
    $$ X \to \prod _{f\colon X\to Y}Y$$
    that takes $x \in X_n$ to $(f(x))_f$, where $f$ ranges over all quotient maps from
$X$ onto nonsingular simplicial sets $Y$.
\end{definition}
\begin{remark}
     Any surjective map from $X$ to a nonsingular simplicial set $Y$ factors uniquely through $DX$. 
\end{remark}
It is not sufficient to only consider the adjoint pair $(D, i)$ when attempting establish a Quillen equivalence from $\sSet$ to $\sSetB$, since the left adjoint $D$ does not induce an equivalence on the homotopy categories; a similar issue arises between the categories $\sSet$ and $\Cat$ \cite[Section 2]{raptisposet}. Desingularization does not preserve homotopy type; for example, desingularization of a loop in $\sSet$ is a point in $\sSetB$. 
A thorough discussion of this homotopical issue, as well as more examples can be found in \cite[Section 4]{fjellbonsSet}.

Subdivision can be used to address this problem, and its interaction with the properties of simplicial sets will be the topic of the remainder of this section. In simplicial sets, subdivision is a generalization of barycentric subdivision of simplicial complexes, in which simplices are partioned into smaller simplices. Note that subdivision of simplicial sets does not change the homotopy type of the geometric realization. We use the following combinatorial definition of simplicial subdivision:

\begin{definition} \cite[Definition 12.2.1]{finitespaces}
\label{sd2}
For a simplicial set $X$, its subdivision, denoted $\Sd X$, consists of the following equivalence classes of tuples as the $q$-simplices:
$$\Sd X_q = (x; S_0, ...S_q)$$
where $x \in X_n$ and $S_i \subset [n]$, and $S_{i} \subset S_{i+1}$ for $0 \leq i < q$. The equivalence relation is specified by:
$$(\mu^*x;S_0,...S_q) \sim (x;\mu(S_0),\mu(S_1)...\mu(S_q)) $$
for a morphism $\mu: [m] \rightarrow [n]$ in $\Delta$, so that $\mu^*x \in X_m$ and $\{S_i\}$ is an increasing list of subsets of $[m]$. $\\$
The simplicial operations are induced by:
$$v^*(x; S_0,...S_q) = (x; S_{v(0)},...S_{v(p)}) $$
for a morphism $v\colon [p] \rightarrow [q]$ in $\Delta$, where $x \in X_n$, and $\{S_i\}$ is an increasing list of subsets of $[n]$. Sd is also a functor between $\sSet$. Given $f\colon X \rightarrow L$ a morphism between two simplicial sets, $f_*= \Sd f\colon \Sd X \rightarrow \Sd L$ is induced by:
$$f_*(x; S_0...S_q) = (f(x); S_0...S_q)$$
\end{definition}

\begin{proposition}
    There is an adjunction
    \[\begin{tikzcd}
            \sSet \arrow[r, shift left=.75ex,"D\Sd^2"{name=G}] & \sSetB\arrow[l, shift left=.75ex, "\Ex^2i"{name=F}] 
            \arrow[phantom, from=F, to=G, "\scriptstyle\dashv" rotate=270] 
        \end{tikzcd}\]
where $i$ is the inclusion, expressing $\sSetB$ as a reflective subcategory of $\sSet$.
\end{proposition}
\begin{proof}
    This is \cite[Lemma 2.2]{fjellbo2020iterative}. 
\end{proof}

\begin{remark}
\label{rmk29}
    Recall that $\Ex$ is the right adjoint to the subdivision functor $\Sd$. In fact, $(\Sd, \Ex)$ is a Quillen equivalence from $\sSet$ to itself (see \cite{dwyerspalinski} and \cite{raptisposet}).
\end{remark}
\subsection{Property C}

 The full subcategory of $\sSetB$ whose objects have Property C is denoted $\sSetBC$. In $\sSetBC$, since each simplex is uniquely determined by its vertices, we will sometimes denote an  $n$-simplex $\sigma \in X_n$ by the ordered list of its vertices $(v_0,v_1,\dots,v_n)$. In Section $\ref{section5}$, we will show that $\sSetBC$ is a reflective subcategory of $\sSetB$.

The following property of $\sSetBC$ is convenient when specifying maps:
\begin{proposition}
\label{vertexmap}
    Any map $f\colon  X \rightarrow Y$ in $\sSetBC$ is uniquely and fully determined by its action on vertices.
\end{proposition}
\begin{proof}

    A map of simplicial sets is fully determined by its action on nondegenerate simplices, by the simplicial relations and the Eilenberg-Zilber Lemma \cite[Theorem 4.2.3]{Fritsch_Piccinini_1990}
    Hence, since a nondegenerate n-simplex in $\sSetBC$ must be the unique nondegenerate simplex with exactly those vertices, specifying $f$ on vertices also gives the map on nondegenerate simplices.
\end{proof}
Another way to state the above property is that the forgetful functor from simplicial sets to the underlying set of vertices is faithful on the subcategory $\sSetBC$.

Definition \ref{sd2} gives a convenient, unique representation of simplices of subdivided simplicial sets:

\begin{definition}[Minimal form of a Subdivided Simplex; \cite{finitespaces}, Definition 12.3.1]
\label{minform}
A $q$-simplex of $\Sd X$ $(x; S_0,...S_q)$ is in minimal form if $x\in X_n$ is nondegenerate and $S_q = [n]$. By Proposition 12.3.2 in \cite{finitespaces}, any simplex in $\Sd X$ can be written in unique minimal form. 
\end{definition}

For future reference when we discuss ordered simplicial complexes, we will need to consider directed loops of 1-simplices which may exist in simplicial sets. In particular, we will see that subdivision results in simplicial sets without any such loops.
\begin{definition}
\label{loopdef}
    An \textit{$n$-loop} is a directed path of $n$ nondegenerate 1-simplices which starts and ends at the same vertex. In particular, a 1-loop has no nondegenerate 1-simplices with the same nondegenerate face.
\end{definition}
\begin{lemma}
\label{propd}
For any simplicial set $X$, $\Sd X$ does not have any $n$-loops.
\end{lemma}
\[
\begin{tikzcd}[cells={nodes={}}]
        \arrow[loop left,distance=3em, start anchor={[yshift=-1ex]west}, end anchor={[yshift=1ex]west}]{}{} \arrow[r, "\Sd"] \bullet 0 \quad
        &\quad  (0,1) \bullet&0 \bullet \quad \arrow[l, bend right,  yshift = -1ex]\arrow[l, bend left, yshift=1ex]
    \end{tikzcd}
    \]
\begin{proof}
    This follows from the minimal form of a nondegenerate simplex in the subdivision. First, no 1-loops exist because a nondegenerate simplex $(x; S_0, S_1)$ requires that $(x; S_0)$ and $(x; S_1)$ can be rewritten as distinct 0-simplices. Now, suppose we had a nondegenerate 1-simplex $(x; S_0, S_1) \in (\Sd X)_1$, with $x \in X_n$ and $S_0$ a $k$-subset of $S_1 = [n]$ specifying a face $y$. Then a 1-simplex in $\Sd X$ which composes with $(x; S_0, S_1)$ can be written in the form $(x'; R_0, R_1)$, where $x' \in X_m$ has $x$ as the $n$-dimensional face specified by $R_0$. Since $x'$ has $x$ as a face, it also has $y$ as a face. Then, we have a 2-simplex $(x'; S'_0, R_0, R_1)$ with $S'_0$ a $k$-subset and $R_0$ an $n$-subset of $R_1 = [m]$. In particular, $d_0(x'; S'_0, R_0, R_1) = (x'; R_0, R_1)$. In the same way, given any nondegenerate $(n-1)$-simplex in minimal form, $(z; S_0...S_{n-1})$, if there is a 1-simplex $(z'; R_0, R_1)$, with $z$ a face of $z'$ specified by $R_0$,  we have an $n$-simplex $(z; S'_0, S'_1..S'_{n-1}, R_0, R_1)$. Hence, any sequence of $n$ composable nondegenerate 1-simplices in $\Sd X$ defines a nondegenerate $n$-simplex, hence it is on the boundary of the simplex. Therefore, the 1-simplices cannot be an $n$-loop.
\end{proof}

\section{Constructing the Left Adjoint of the Inclusion of $\sSetBC$ into $\sSetB$}
\label{section5}
A simplicial set in $\sSetB$ may fail to have Property C if it has two 1-simplices with vertices $\{v_0, v_1\}$ but with opposing source and target, or two parallel $n$-simplices with shared boundary. This section establishes that $\sSetBC$ can be thought of as the localization of $\sSetB$ where opposing 1-simplices are collapsed to a point and parallel simplices to a single simplex. 

With this in mind, the main result of this section, Theorem \ref{mainthm}, establishes an adjunction between $\sSetB$ and $\sSetBC$ as a reflective localization, by introducing a left adjoint to the fully faithful inclusion functor $i$. An illustration of $L$ for some simple cases is given below.
\[
    \includegraphics[width=1\linewidth]{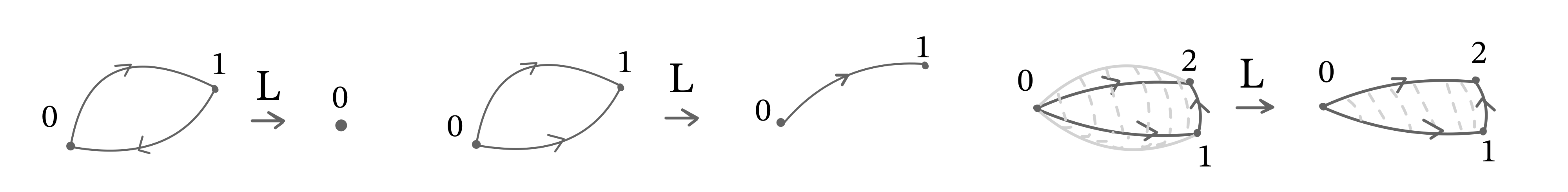}
\]
   On left is the simplicial set consisting of two opposing nondegenerate 1-simplices, which is collapsed to a point. The middle and right are the cases of parallel nondegenerate 1 and 2-simplices, which are collapsed to a single copy of the respective simplex.

Our construction will use the small object argument in a similar way to the proof of the factorization axiom for model categories established in \cite[Section 7.12]{dwyerspalinski}.

\begin{notation}
   We define a family of simplicial set maps $\mathcal{F} = \{f_i\colon  A_i \rightarrow \Delta[i]\}$ with $A_i$ defined as follows. Let $A_0 $ be the simplicial set of opposing 1-simplices, $0 \to 1 \rightarrow 0$, and let $f_0: A_0 \rightarrow 0$ be the collapse of 2 opposing 1-simplices. For $n > 0$, let $A_n$ be the simplicial set $\Delta[n] \cup_{\partial\Delta[n]} \Delta[n]$, where the union is over the (shared) boundary. Let $f_n\colon  A_n \rightarrow \Delta[n]$ be the quotient to the $n$-simplex. 
   \label{defmaps}
\end{notation}
The family $\mathcal{F}$ of Notation \ref{defmaps} can be used to identify the image of the inclusion $i\colon \sSetBC \hookrightarrow \sSetB$ via a right lifting property, which we now show.
\begin{proposition}
A simplicial set $X$ with Property B has Property C if and only if $X \rightarrow *$ has the right lifting property relative to each of the maps in $\mathcal{F}$. 
\label{lemma123}
\end{proposition}
\begin{proof}

First suppose $X$ has Properties B and C; we would like to show $X \rightarrow *$ has the right lifting property with respect to all maps $f_i \in \mathcal{F}$. 
\[ \begin{tikzcd}
 A_i \arrow{r}{g} \arrow[swap]{d}{f_i} & X \arrow{d}{} \\%
 B_i \arrow[ur,"\exists h",dashed]\arrow{r}{}& *
\end{tikzcd}
\]

Given a diagram as above, we will construct a map $h$ such that $g=hf_i$. Starting with the case $n=0$, recall that $f_0$ is the constant map $(0 \rightarrow 1 \rightarrow 0)\to 0$.
As $g$ is a map of simplicial sets, we must have 1-simplices $g(0) \rightarrow g(1)$ and $g(1) \rightarrow g(0)$ in $X$. However, since $X$ has Property B and C, these simplices must be degenerate
(since otherwise we would have two nondegenerate 1-simplices with the same vertices, a contradiction), forcing $g(0) = g(1)$. We can thus send the single element in $B_0$ to this simplex so that $hf_0(A_0) = g(A_0)$. 

Similarly, for $f_n\colon  \Delta[n] \cup_{\partial\Delta[n]} \Delta[n] \rightarrow \Delta[n]$, we notice that the $n$-simplices provided by each copy of $\Delta[n]$ share a boundary, and hence share vertices. Suppose $\sigma_1, \sigma_2$ are these two $n$-simplices sharing the distinct vertices $(v_1, ...v_{n+1})$. Since $g$ is a simplicial map, the $i$th vertex of $g(\sigma_1)$ and $g(\sigma_2)$ are the same, and hence their images have the same tuple of $n$-vertices. Then these simplices have the same underlying nondegenerate simplex $\sigma$ consisting of the same distinct vertices. The Eilenberg-Zilber Lemma states that any simplex of a simplicial set can be uniquely written as $s^{\sharp}\sigma$ for a degeneracy operator $s^\sharp$, and nondegenerate simplex $\sigma$, so it suffices to show that $s^\sharp$ is the same. Since the tuples of vertices are identical, this follows directly from the simplicial relations $s_is_j = s_js_{i-1}$ for $i > j$, (which also implies $s_i^2 = s_{i+1}s_i$). Hence we can equate any degenerate simplex by duplicating particular vertices to obtain the same $n+1$-tuple, so that any two $n$-simplices in a simplicial set with Properties B and C with the same $(n+1)$-tuple of vertices are the same simplex. We can then send $\Delta[n]$ in $B_n$ to the uniquely determined image $g(A_n)$ so that $hf_n(A_n)=g(A_n)$.


Conversely, we would like to show that if $X \rightarrow *$ has the right lifting property relative to each $f_i \in \mathcal{F}$, then $X$ has Property C, provided that it has Property B.

Note that the fact that $X\to *$ has the right lifting property relative to $f_0$ implies that $X$ has no opposing nondegenerate 1-simplices; similarly, its lifting property relative to $f_1$ implies that $X$ has no parallel 1-simplices. Hence $X$ satisfies  Property C on all 1-simplices.

Suppose $X$ satisfies Property C on all $(n-1)$-simplices, and let $\sigma_1, \sigma_2\in X_n$ be two nondegenerate simplices with the same vertices. Then $\sigma_1$ and $\sigma_2$ 
must have the same boundary; this is since $X$ is in $\sSetB$, all nondegenerate simplices are nonsingular, and its $n-1$ dimension simplices are all uniquely determined by satisfying Property C. Then, 
this determines a simplicial map $\Delta[n] \cup_{\partial\Delta[n]} \Delta[n] \rightarrow X$ mapping to these simplices, and the lifting property then implies that $\sigma_1=\sigma_2$ in $X$.
\end{proof}
Now we can introduce the left adjoint to the inclusion $i\colon  \sSetBC \hookrightarrow \sSetB$.
\begin{construction}
We define a functor $L\colon  \sSetB \to \sSetBC$ as follows. Given an object $X \in \sSetB$, we can use the family $\mathcal{F}$ of Notation 
\ref{defmaps} and the small object argument to factor $X \to *$ as $X \to L(X) \to *$ so that $L(X) \to *$ has the right lifting property relative to $\mathcal{F}$. By Proposition \ref{lemma123}, this implies that $L(X) \in \sSetBC$. Furthermore, the factorization defining $L$ is functorial, since this is the weak factorization system cofibrantly generated by $\mathcal{F}$, which is functorial (see Remark 2.3 from \cite{guetta2023fibrantlyinducedmodelstructures} and Proposition 16 of \cite{Bourke_2016}). 
\label{deffunctor}
\end{construction}

\begin{theorem}
Let $L$ be as defined in Construction \ref{deffunctor}. Then $\sSetB\colon  L \dashv i : \sSetBC$ is an adjunction and hence $\sSetBC$ is a reflective subcategory of $\sSetB$.
\label{mainthm}
\end{theorem}
\begin{proof}
To show $L$ is left adjoint to the inclusion $i\colon  \sSetBC \hookrightarrow \sSetB$, consider the natural transformation $\eta: 1_{\sSetB} \rightarrow iL $ with component $\eta_X: X \to iLX$; this component is an element of Hom$_{\sSetB}(X,iLX)$. As stated in Construction \ref{deffunctor}, we have the naturality of the following diagram for all $X \xrightarrow{f} X'$ in $\sSetB$:
\[ \begin{tikzcd}
X \arrow{r}{\eta_X} \arrow[swap]{d}{f} & iLX \arrow{d}{iL(f)} \\%
X' \arrow{r}{\eta_{X'}}& iLX'
\end{tikzcd}
\]
Now we define the natural transformation $\varepsilon\colon  Li \rightarrow 1_{\sSetBC}$, for $Y \in \sSetBC$, with component $\varepsilon_Y$ equal to $1_Y$. Since $1_Y$ is an isomorphism for each $Y$, this means $\varepsilon$ is a natural isomorphism, so that $Li$ is naturally isomorphic to the identity.

Now, we need to verify the triangle identities. As established, $Li$ is isomorphic to the identity by definition. Then the triangle identities reduce to:
\[
\begin{tikzcd}[row sep=3em]
LX \arrow{r}{L\eta_X} \arrow[""{name=foo}]{dr}[swap]{1_{LX}} & LX \arrow{d}{\varepsilon LX} \arrow[Rightarrow, from=foo, swap, near start, "\cong"]& iY \arrow{d}[swap]{\eta iY} \arrow[""{name=bar, below}]{dr}{1_{iY}} & \\
& LX & iY \arrow{r}[swap]{i\varepsilon_Y} \arrow[Rightarrow, to=bar, swap, near start, "\cong"] & iY
\end{tikzcd}\] 
It is almost immediate that the triangle legs are identities, hence the compositions are both isomorphic to $1_{LX}$ and $1_{iY}$ respectively.
\end{proof}

\section{Right Transferring the Model Structure from $\sSetB$ to $\sSetBC$}
In this section we establish a cofibrantly generated model structure on $\sSetBC$; we do this by lifting the model structure from $\sSetB$, due to \cite{fjellbonsSet}, along the adjunction established in Section 3, Theorem \ref{mainthm}, to $\sSetBC$. To do so, we recall some background for cofibrantly generated model categories.

\begin{definition}
\label{relcell}
Let $I \subset Mor(\mathcal{C})$ be a set of morphisms in a category $\mathcal{C}$. Denote $cell(I)$ as the relative $I$-cell complexes, the class of morphisms in $\mathcal{C}$ obtained by transfinite composition of pushouts of elements in $I$.

\end{definition}

We will use the following result for right-transfer of cofibrantly-generated model categories:

\begin{theorem}[9.1: Lifting Lemma in \cite{DwyerHirschhornKanSmith1997}]
\label{thm1}
Let $\mathcal{M}$ be a cofibrantly generated model category with generating cofibrations $I$ and generating trivial cofibrations $J$. Let $\mathcal{N}$ be a category that is closed under small limits and colimits, and let \[
        \begin{tikzcd}
            \mathcal{M} \arrow[r, shift left=.75ex,"F"{name=G}] & \mathcal{N}\arrow[l, shift left=.75ex, "U"{name=F}] 
            \arrow[phantom, from=F, to=G, "\scriptstyle\dashv" rotate=270]
        \end{tikzcd}
    \] be a pair of adjoint functors. If
\begin{enumerate}
    \item Both of the sets FI and FJ permit the small object argument.
    \item U takes cell(FJ) to weak equivalences.
\end{enumerate} 
Then there is a cofibrantly generated model category structure on $\mathcal{N}$ in which FI is a set of generating cofibrations, FJ is a set of generating trivial cofibrations, and the weak equivalences are the maps that U takes into a weak equivalence in $\mathcal{M}$. Furthermore, with respect to this model category structure, (F, U) is a Quillen pair.
\end{theorem}
\begin{notation}
\label{def85}
    For the rest of this section we refer to $I = \{\partial\Delta[n] \hookrightarrow \Delta[n]| n \geq 0\}$ and $J = \{\Lambda^k[n] \hookrightarrow \Delta[n]|n \geq 0, 0 \leq k \leq n\}$, which are the generating cofibrations and trivial cofibrations of $\sSet$. 
\end{notation}

The category $\sSetB$ itself has a cofibrantly generated model structure which is Quillen equivalent to the standard model structure on $\sSet$, which we restate below:

\begin{theorem}[\cite{fjellbonsSet}, Theorem 1.2]
    There is a proper, cofibrantly generated model structure on $\sSetB$ right-transferred along the adjunction
     \[
        \begin{tikzcd}
           \sSet \arrow[r, shift left=.75ex,"D\Sd^2"{name=G}] & \sSetB\arrow[l, shift left=.75ex, "\Ex ^2i"{name=F}] 
            \arrow[phantom, from=F, to=G, "\scriptstyle\dashv" rotate=270]   
        \end{tikzcd}
    \]
    such that the Quillen pair $(D\Sd^2, \Ex^2i)$ is a Quillen equivalence.
    \label{thm122}
\end{theorem}

By Theorem \ref{thm122}, there is a cofibrantly generated model structure on $\sSetB$ which has generating cofibrations $D\Sd^2I$ and generating trivial cofibrations $D\Sd^2J$. We apply Theorem $\ref{thm1}$ to $\sSetB$ to lift the model structure of Theorem \label{thm12} to $\sSetBC$ along the adjoint pair $(L,i)$ of Theorem \ref{mainthm}. For this, it helps to understand first which objects are small in these categories.

\begin{lemma}
An object in $\sSetB$ or $\sSetBC$ is small if and only if it has a finite number of non-degenerate simplices.
    \label{slemma}
\end{lemma}
\begin{proof}
Given an object $X\in\sSetB$, we want to show that $X$ has a finite number of nondegenerate simplices precisely when, for any filtered diagram $A\colon \mathcal{J}\to\sSetB$, we have:
$$\Hom_\sSetB(X, \colim_j(A_j))\cong \colim_j(\Hom_\sSetB(X,A_j))$$
Let $X\in\sSetB$ have a finite number of nondegenerate simplices. Then $X$ is a small object in $\sSet$, and so it suffices to show that $\Hom_\sSet(X, \colim_j(iA_j))\cong \Hom_\sSetB(X, \colim_j(A_j))$ whenever $A\colon \mathcal{J}\to \sSetB$ is a filtered diagram.

Suppose otherwise, that we had some $n$-simplex $\sigma$ with repeated vertices in $\colim_j (iA_j)$. Consider the map $\sigma\colon\Delta[n] \to \colim_j (i A_j)$. Note that $\Delta[n]$ is small in $\sSetB$ and $\sSetBC$, since any simplex in the colimit must factor through some $A_j$. Hence, $\sigma$ factors through some $A_j$, and so we must have some $k > j$ and $A_j \to A_{k}$ which sends a nonsingular simplex to a nondegenerate singular one; but then not all $A_j$ are in $\sSetB$, which contradicts the fact that all $A_j$ are in $\sSetB$. For $\sSetBC$ we have the same argument, but take $\Delta[n]$ to be disjoint union of $\Delta[n]$ instead, with the map identifying their boundaries in the colimit. A similar argument holds for $\sSetBC$, replacing the map $\Delta[n]\to\colim_j (i A_j)$ by a map $\Delta[n]\sqcup\Delta[n]\to \colim_j (i A_j)$ that identifies the boundaries of the two $n$-simplices.

 The reverse direction follows the same argument as compactness of finite sets. We repeat the argument here for $\sSetB$: any object $X\in\sSetB$ can be expressed as a filtered colimit of its finite subsimplicial sets $X_i$ (i.e.\ ordered by inclusion with a finite number of nondegenerate simplices), so we can take $X_j$ as the filtered colimit of subsimplices, starting from a finite set of simplices in $X_0$. Each of these subsimplices $X_j$ are in $\sSetB$ since $X \in \sSetB$. Then, since $D$ preserves colimits, this is still the colimit in $\sSetB$, as $DX= X \in \sSetB$. If $X$ is small, then the map $\id_X\colon X\to \colim_j X_j$ must factor as $X\to X_k\to \colim_j X_j$ for some $k$, which implies that $X$ has a finite number of nondegenerate simplices, as $X_k$ does. 
\end{proof}
 
The following definition and lemmas will be needed to establish a model structure for $\sSetBC$, as well as when passing to ordered simplicial complexes. 
\begin{definition}
   An inclusion of simplicial sets $j\colon  X \to Y$ is a \textit{full simplicial inclusion} if any simplex of $Y$ whose vertices are in $Im(j)$ is a simplex in $Im(j)$.
\end{definition}
\begin{remark}
    Although the following result and proof follows for full simplicial inclusions in complete generality, we will apply it to actual inclusions of simplicial sets. As such, we will abuse notation when convenient and treat full simplicial inclusions as honest inclusions.
\end{remark}
 \begin{lemma}
      Pushouts preserve full simplicial inclusions in $\sSetB$.
     \label{lemma2}
 \end{lemma}
 \begin{proof}
 Suppose we had the following pushout diagram on the left in $\sSetB$, where $j$ is a full simplicial inclusion:
\[ \begin{tikzcd} X\arrow[r, "f"] \arrow[swap]{d}{j} & Z \arrow[d, dashed, "p^n"] \\%
 Y \arrow[r, dashed, "g^n"]& P^n
\end{tikzcd} \quad \quad\quad\begin{tikzcd} X\arrow[r, "f"] \arrow[swap]{d}{j} & Z \arrow[d, dashed, "p^s"] \\%
 Y \arrow[r, dashed, "g^s"]& P^s
\end{tikzcd}  
\]     
To obtain the pushout in $\sSetB$, since $\sSetB$ is a reflective subcategory of $\sSet$ \cite{fjellbonsSet}, it suffices to first find the pushout $P^s$ in $\sSet$, pictured on the right, and then desingularize. 

We show that $p^s: Z \to P^s$ is a full simplicial inclusion. Recall that $P^s$ can be computed as the simplicial set whose $n$-simplices are given by $P^s_n = (Z_n\sqcup Y_n)/f(x) \sim j(x)$ for all $x \in X_n$. Suppose there is a simplex $\sigma\in P^s$ with vertices $\{p^s(v_0),\dots,p^s(v_n)\}$, where $\{v_0, \dots v_n\}$ are  vertices in $Z$. By definition of $P^s$, $\sigma$ must be the image of an $n$-simplex from $Y$ or $Z$; if $\sigma$ is the image of an $n$-simplex from $Z$, then we are done. Otherwise, $\sigma = g^s(\sigma')$ is either nondegenerate simplex or the degeneracy of a nondegenerate simplex $\bar{\sigma}$, expressed as $\sigma = s\bar{\sigma}$, and we can apply the following argument to $\bar{\sigma}$. Suppose first that $\sigma$ is a nondegenerate simplex, so that $\sigma'$ is also nondegenerate. Since $Y \in \sSetB$, $\sigma'$ has distinct vertices $\{v'_0, \dots v'_n\}$ where $g^s(v_i') = p^s(v_i)$. Since $g^s(v_i') \in Y$ are identified with $p^s(v_i)$ in $P^s$, there must exist distinct vertices $\{v'_0, \dots v'_n\}$ in $X$. Since $j$ is a full simplicial inclusion, there exists $\hat{\sigma} \in X$ such that $j(\hat{\sigma}) = \sigma'$. As an inclusion, $j$ is a monomorphism, and hence its pushout $p^s$ in $\sSet$ is a monomorphism; therefore, the vertices of $f(\hat{\sigma})$ must be $\{v_0,\dots,v_n\}$ since they are mapped by $p^s$ to $\{p^s(v_0),\dots,p^s(v_n)\}$, so that $p^s: Z \to P^s$ is a full simplicial inclusion in $\sSet$. 

Now, we want to show that $p^n= Dp^s: Z \to P^n = DP^s$ is a full simplicial inclusion in $\sSetB$. Since $D$ only identifies simplices, $p^s$ is full,  and therefore $DP^s$ has fewer simplices than $P^n$. It is an inclusion because $p^s$ factors as $Z \xrightarrow{Dp^s}DP^s \xrightarrow{i}P^s$ (since a nondegenerate simplex in $Z$ has distinct vertices so that $p^s$ as an inclusion maps this to a simplex with distinct vertices), where both $i$ and $p^s$ are monomorphisms. This implies that $Dp^s$ is also, and is hence an inclusion.
 \end{proof}

 \begin{lemma}
      Pushouts preserve full simplicial inclusions in $\sSetBC$.
     \label{lemma3}
 \end{lemma}
 \begin{proof}
     In Lemma $\ref{lemma2}$, we have established that pushouts in $\sSetB$ are preserved along full simplicial inclusions. Starting from a span as in the preceding lemma, consider the pushout in $\sSetB$, $P^n= \colim(Y \xleftarrow{j} X \xrightarrow{f} Z)$ where $X, Y, Z$ are in $\sSetBC$. Pushouts are $P = LP^n$ in $\sSetBC$ since they are obtained by first taking pushouts of the diagram inclusions in $\sSetB$ and then applying $L$,
     
      Applying $L$ either does not affect or reduces the number of simplices via a quotient relation identifying parallel simplices or identifying opposing edges to a single vertex. It is immediate that the pushout $p\colon  Z \to LP^n$ is still full; to see that it is still an inclusion, note that $p^n$ factors as $Z \xrightarrow{Lp^n}LP^n \xrightarrow{i}P^n$ , where both $i$ and $p^n$ are monomorphisms. This factorization exists because any simplex in $P^n$ with vertices in $Z$ must be mapped to by the \textit{unique} simplex in $Z$ with those vertices, due to $Z$ being in $\sSetBC$ and $p^n$ a full simplicial inclusion. This implies that $Lp^n = p$ is also a monomorphism, and is hence an inclusion.
    
 \end{proof}

 To verify the second condition of Theorem $\ref{thm1}$, which we show in Lemma $\ref{mainlemma}$, we need the following lemmas which construct a deformation retraction will be used to show the required weak equivalence.
\begin{lemma}
     Let $m \geq 2$. There is a deformation retraction $h$ from $|\Sd^m\Delta[n]|$ to $|\Sd^m\Lambda^k[n]|$, such that given a quotient relation $\sim$ on $\Sd^m\Lambda^k[n]\subset \Sd^m\Delta[n]$, there is an associated deformation retraction $h': |D(\Sd^m\Delta[n]/\sim)| \to |D(\Sd^m\Lambda^k[n]/\sim)|$.
     \label{defret}
\end{lemma}
\begin{proof}

    We show the statement for $m =2$ and the same argument immediately generalizes to all $m \geq2$. Let $x\in |\Sd^2\Delta[n]|$; we can consider the geometric realization of the (subdivided) standard simplex and the horn to be embedded in Euclidean space. We take the deformation retraction as following the linear path formed by projection to one of the coordinates, starting from $x$. Explicitly, in $\mathbb{R}^m$, let $p_x$ be the point of intersection of the parallel line to $|\Sd^2\Lambda^k[n]|$ by projecting to the $k$th horn, and the deformation retraction from $|\Sd^2\Delta[n]|$ to $|\Sd^2\Lambda^k[n]|$ is given by $h_t(x)= (p_x-x)t+x$, for $t \in [0,1]$.  For an example on $\Delta[2]$ see Figure \ref{89}.

    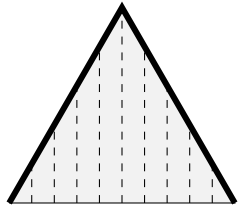
\begin{figure}[h]
    \begin{tikzpicture}[scale=0.7]
          \coordinate (A) at (0,0);
    \coordinate (B) at (3, 0);
    \coordinate  (C) at (1.5,2.6);
   
    \fill [gray!10](A) --(B)--(C) --cycle;
    \draw[line width=0.8mm] (A) --(C)--(B);
    \draw[dashed] (0.3, 0) --(0.3,.52);
    \draw[dashed] (0.6, 0) --(0.6,1.04);
    \draw[dashed] (0.9, 0) --(0.9,1.56);
    \draw[dashed] (1.2, 0) --(1.2,2.08);
    \draw[dashed] (1.5, 0) --(1.5,2.6);
    \draw[dashed] (1.8, 0) --(1.8,2.08);
    \draw[dashed] (2.1, 0) --(2.1,1.56);
    \draw[dashed] (2.4, 0) --(2.4,1.04);
    \draw[dashed] (2.7, 0) --(2.7,.52);
    \draw(A) --(B)--(C);
    \end{tikzpicture}
  

    \caption{The deformation retraction obtained by projecting to the $k$th horn for $|\Delta[2]|$. The dashed lines indicate the parallel lines along which the $p_x$ are defined.}
    \label{89}
    \end{figure}

      By the second subdivision, any nondegenerate simplex in $\Sd^2\Delta[n]$ is the face of an $n$-simplex or is itself one, so we can restrict our attention to these. The quotient boundary relation $\sim$ can identify simplices of $\Sd^2\Lambda^k[n]$ only, so no vertices in  $\Sd^2\Delta[n]\backslash\Sd^2\Lambda^k[n]$ are identified. Thus, $n$-simplices in $\Sd^2\Delta[n]$ without multiple vertices on $\Sd^2\Lambda^k[n]$ must remain nondegenerate after quotienting; call these simplices unaffected. An illustration of an example of these simplices is shown for $n = 2$ in Figure 2. Note that the unaffected simplices form a connected component in $\Sd^2\Delta[n]$, and that due to the second subdivision, any vertex on $\Sd^2\Lambda^k[n]$ is within one 1-simplex away from the vertex of an unaffected simplex.

   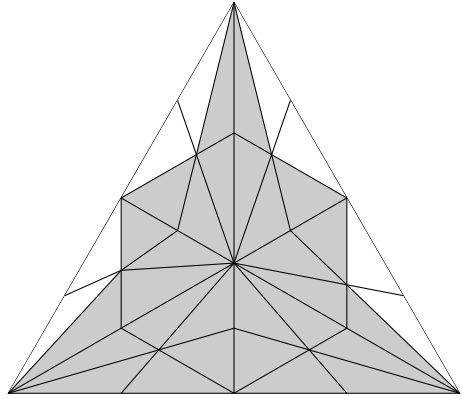
\begin{figure}[h]
\begin{tikzpicture}[scale=0.5]
 \coordinate (A) at (0,0);
    \coordinate (B) at (6, 0);
    \coordinate  (C) at (3,5.2);
    \fill[gray!40] (A)--(B)--(C)--cycle;
     \draw (A) -- (B) -- (C) -- cycle;
     \coordinate (F) at (3,0); 
    \coordinate (G) at (1.5, 2.598); 
    \coordinate  (H) at (4.5,2.598);   
    \coordinate  (I) at (3, 1.732);
     \draw (A) -- (I);
     \draw (B) -- (I);
     \draw (C) -- (I);
     \draw (F) -- (I);
     \draw (G) -- (I);
     \draw (H) -- (I);
     \coordinate (A1) at (1.5,0.866);
    \coordinate (B1) at (1.5, 0);
    \coordinate  (C1) at (3,0.866); 
    \coordinate (M1) at (2, 0.577);
     \draw (A1) -- (M1);
     \draw (B1) -- (M1);
     \draw (C1) -- (M1);
     \draw (F) -- (M1);
     \draw (I) -- (M1);
     \draw (A) -- (M1);
     \coordinate (A2) at (4.5,0.866);
    \coordinate (B2) at (4.5, 0);
    \coordinate  (C2) at (3,0.866); 
    \coordinate (M2) at (4, 0.577);
     \draw (A2) -- (M2);
     \draw (B2) -- (M2);
     \draw (C2) -- (M2);
     \draw (F) -- (M2);
     \draw (I) -- (M2);
     \draw (B) -- (M2);
     \coordinate (A3) at (1.5,0.866);
    \coordinate (B3) at (0.750,1.299);
    \coordinate  (C3) at (2.250,2.165); 
    \coordinate (M3) at (1.5, 1.6357);
            \fill[white] (M3) -- (A)--(G) ;
     \draw (A3) -- (M3);
     \draw (B3) -- (M3);
     \draw (C3) -- (M3);
     \draw (A) -- (M3);
     \draw (I) -- (M3);
     \draw (G) -- (M3);
     \coordinate (A4) at (3.000,3.464);
    \coordinate (B4) at (2.25,3.897);
    \coordinate  (C4) at (2.250,2.165); 
    \coordinate (M4) at (2.5, 3.175);
            \fill[white] (M4) -- (C)--(G) ;
     \draw (A4) -- (M4);
     \draw (B4) -- (M4);
     \draw (C4) -- (M4);
     \draw (C) -- (M4);
     \draw (I) -- (M4);
     \draw (G) -- (M4);
     \coordinate (A5) at (4.500,0.866);
    \coordinate (B5) at (5.250,1.299);
    \coordinate  (C5) at (3.750,2.165); 
    \coordinate (M5) at (4.5, 1.4433);
        \fill[white] (M5) -- (B)--(H) ;
     \draw (A5) -- (M5);
     \draw (B5) -- (M5);
     \draw (C5) -- (M5);
     \draw (B) -- (M5);
     \draw (I) -- (M5);
     \draw (H) -- (M5);

     \coordinate (A6) at (3,3.464);
    \coordinate (B6) at (3.750, 3.897);
    \coordinate  (C6) at (3.750,2.165); 
    \coordinate (M6) at (3.5,3.1753 );
    \fill[white] (M6) -- (C)--(H) ;
     \draw (A6) -- (M6);
     \draw (B6) -- (M6);
     \draw (C6) -- (M6);
     \draw (C) -- (M6);
     \draw (I) -- (M6);
     \draw (H) -- (M6);
     
\end{tikzpicture}
        \caption{An illustration for $\Sd^2\Delta[2]$ where the horn excludes the bottom edge; the shaded simplices denote simplices which always stay nondegenerate after identification. Notice that these simplices specify a path-connected component in the interior of $\Sd^2\Delta[2].$}
        \label{22}  
    \end{figure}
      
      Now consider the second geometric barycentric subdivision of $\Delta^n$, where to each subdivided geometric $n$-simplex we can associate the corresponding $n$-simplex of $\Sd^2\Delta[n]$. Let $X'$ be the subsimplex consisting of the connected component of unaffected simplices.  Note that the $n$-simplices in $\Sd^2\Delta[n]$, but not fully contained in $X'$, completely cover $\Sd^2\Lambda^k[n]$, and furthermore the collection of lines defined from $p_x$ must intersect one such $n$-simplex. Thus after passing to the quotient, the collection of parallel lines still maps onto $|\Sd^2\Lambda^k[n]/\sim|$, and we can define a deformation retraction $h'$ from the quotient spaces. Let $x \in |\Sd^2\Delta[n]|$ and $[x]$ its corresponding point in $|D(\Sd^2\Delta[n]/\sim)|$. We define $h'$ for each $[x] \in |D(\Sd^2\Delta[n]/\sim)|$ by $h'([x], t) = [h(x, t)]$.
Note that $[x]$ is in one-to-one correspondence with $x$ away from identified simplices (where the identified values of $x$ are collapsed along the vertical paths), so $h'$ is well-defined.

      To see that $h'$ is in fact a deformation retraction onto $|\Sd^2\Lambda^k[n]/\sim|$, note that any point $p$ in $|\Sd^2\Lambda^k[n]/\sim|$ is either in the interior of the geometric realization of a nondegenerate $(n-1)$-simplex in $|\Sd^2\Lambda^k[n]/\sim|$, or it is on the boundary of such a nondegenerate simplex. There is a sufficiently small neighborhood of $p$ contained in the geometric realization of one or more $n$-simplices in $\Sd^2\Delta[n]$ (there may be multiple if $p$ is on a boundary). In either case, we may select a point $\tilde{x}$ in the geometric realization of one of these nondegenerate simplices such that the point $p_{\tilde{x}}$ coincides with $p$.

\end{proof}

\begin{lemma}
     Given the deformation retraction $h'$ from the statement of Lemma \ref{defret} with $m \geq 2$, there is an associated deformation retraction $h'': |LD(\Sd^m\Delta[n]/\sim)| \to |LD(\Sd^m\Lambda^k[n]/\sim)|$.
     \label{defret2}
\end{lemma}
\begin{proof}
We again show the statement for $m =2$, and the same argument follows for all $m \geq 2$. Consider $h'$, the deformation retraction from $| D(\Sd^2\Delta[n]/\sim)|$ to $| D(\Sd^2\Lambda^k[n]/\sim)|$ as defined in Lemma \ref{defret}. We will see that $h''$ can be defined similarly, in that we only need to consider simplices adjacent to the boundary.

 If an $n$-simplex (of the second subdivided $n$-simplex) only has one vertex on $\Sd^2\Lambda^k[n]$, then due to the second subdivision it is the unique $n$-simplex with the same $n$ vertices in $\Sd^2\Delta[n]\backslash\Sd^2\Lambda^k[n]$, such that the final vertex is on $\Sd^2\Lambda^k[n]$. 
 
 This implies that since $\sim$ and $D$ only identify vertices on the boundary, that we only need to consider simplices with multiple vertices on $\Sd^2\Lambda^k[n]$, since $L$ identifies simplices with all the same vertices.


These simplices are precisely those that may have vertices identified, i.e. simplices that may become degenerate, after applying $D$ (the simplices that are not considered unaffected). Note that these share faces with simplices that must stay nondegenerate after identification.

In light of this, we can define $h''$ in essentially the same was as $h'$: take $x \in |\Sd^2\Delta[n]|$, and define $h''([x], t)$ along the parallel line defined in Lemma \ref{defret}. To see this, note that if $x\in\sigma$ where $\sigma$ becomes identified after applying $LD$, then either $\sigma$ becomes the degeneracy of a nondegenerate $m$-simplex, a face of another nondegenerate $n$-simplex in $\Sd^2\Delta[n]$, or $\sigma$ is identified with another $n$-simplex with the same boundary after applying $L$ (if it is identified by $L$). In the latter case, we can simply choose one of the parallel simplices that $\sigma$ is identified with, since the boundary is the same. Then as before, we can define $h''([x], t) = [h(x, t)]$. 
\end{proof}
\begin{lemma}
\label{thmbc}
   For $n \geq 0$ and $0 \leq k \leq n$, the simplicial sets $D\Sd^2\Lambda^k[n]$ and $D\Sd^2\Delta[n]$ in $\sSetB$ have Property C.
\end{lemma}
\begin{proof}
    The simplicial sets $\Delta[n]$ and $\Lambda^k[n]$ have Property B and thus every face of a nondegenerate simplex is nondegenerate \cite[Theorem 12.1.5]{finitespaces}. 
    Thus by \cite[Theorem 12.4.1]{finitespaces}, $\Sd^2\Delta[n]$ and $\Sd^2\Lambda^k[n]$ have Properties B and C. Furthermore, the desingularization functor $D$ does not change $\Sd^2\Delta[n]$ or $\Sd^2\Lambda^k[n]$ since they already have Property B.
\end{proof}
\begin{remark}
    We refer to the generating trivial cofibrations, $\{D\Sd^2\Lambda^k[n] \hookrightarrow D\Sd^2\Delta[n]\mid n \geq 0, 0 \leq k\leq n$\}, as $\{\Sd^2\Lambda^k[n] \hookrightarrow \Sd^2\Delta[n]\mid n \geq 0, 0 \leq k\leq n$\} since the latter are already completely contained in $\sSetB$ by Lemma $\ref{thmbc}$. 
\end{remark}
 \begin{lemma}
    Let $m \geq 2$, and $j' \in \{ \Sd^m\Lambda^k[n] \hookrightarrow \Sd^m\Delta[n]\mid n \geq 0, 0 \leq k\leq n\}$ be a subdivision of the generating trivial cofibrations of $\sSetB$. Let $f\colon  \Sd^m\Lambda^k[n] \to Y$ be a map in $\sSetBC$.
    Then the pushout $p$ in $\sSetBC$:
     \[
\begin{tikzcd}
\Sd^m\Lambda^k[n] \arrow[d,"j'"'] \arrow[r, "f"]&  Y \arrow[d, "p"']  \\
\Sd^m\Delta[n] \arrow[r]& P 
\end{tikzcd}
\]
is mapped by the inclusion to a weak equivalence in $\sSetB$. 
\label{mainlemma}
\end{lemma}
\begin{proof}
For notational convenience, we let $j': X \to Z$, Furthermore, $X, Y$, and $Z$ will be referred to as such in both $\sSetBC$ and $\sSetB$; this is because $iLX =X, iLZ = Z$, as $X$ and $Z$ already have Properties B and C (so we can identify them as an object in $\sSetBC$).

Denote the pushout of the underlying diagram in $\sSetB$ as $P' := \colim(Z \xleftarrow{j'} X \xrightarrow{f}Y)$, the map to the pushout $p'\colon Y \to P'$ in $\sSetB$. 

Since $j'$ is a generating trivial cofibration in $\sSetB$, $p'$ is a weak equivalence in $\sSetB$. The pushout $p\colon  Y \to P$ in $\sSetBC$ is obtained by taking $p'\colon Y \to P'$ and then applying $L$, since $\sSetBC$ is a reflective subcategory of $\sSetB$. 

To first compute the pushout $p'\colon Y \to P'$ in $\sSetB$, note that we can take the pushout of the underlying diagram in $\sSet$ and apply $D$ as $\sSetB$ is a reflective subcategory of $\sSet$.  The pushout in $\sSet$ can be viewed as attaching $Z$ to $Y$ along $f\colon  X \to Y$; this is a weak equivalence $\phi\colon  Y \to Z \sqcup_f Y$ (since $D\phi$ is a weak equivalence in $\sSetB$). Then $p = Lp' = LD\phi\colon  L DY \to L D(Z \sqcup_f Y)$.  By the induced model structure on $\sSetB$, a map in $\sSetB$ is a weak equivalence if and only if it is a weak equivalence in $\sSet$, if and only if it is a weak homotopy equivalence in the geometric realization. Hence it suffices to show that $p = LD\phi$ induces a homotopy equivalence between $|LDY| \cong |Y|$ (since $Y \in \sSetBC$) and $|L D(Z \sqcup_f Y)| $. Since full simplicial inclusions are preserved under pushouts in both $\sSetB$ and $\sSetBC$ by Lemmas \ref{lemma2} and \ref{lemma3}, $Y \hookrightarrow LD(Z \sqcup_f Y)$ is a full simplicial inclusion. Thus the map $|LD\phi|$ is the inclusion of $|Y|$ into $|LD(Z \sqcup_f Y)|$, and so if $|\tilde{\phi}|$ is a retraction onto $|Y|$, then $|\tilde{\phi}||LD\phi| = \id_{|Y|}$. Then, to show $|LD\phi||\tilde{\phi}| \simeq \id_{|LD(Z \sqcup_f Y)|}$, it suffices to give an appropriate deformation retraction of $|LD(Z \sqcup_f Y)|$ onto $|Y|$.
Then, it suffices to show that if $\sim$ is the equivalence relation on $X$ induced by the pushout in $\sSet$, then $|LD(Z/\sim)|$ deformation retracts onto $|LD(X/\sim)| \cong |\text{Im}(f)|\subset |Y|$. This is because in order for a simplex to be identified by $D$ or $L$, it must be provided from $Z$ since $P$ must have all the simplices from $Y$, as $p$ is a full simplicial inclusion by Lemma \ref{lemma3}. 
By Lemma \ref{defret2}, the result follows, so that $|LD(Z \sqcup_f Y)|$ deformation retracts onto $|Y|$, and $LD\phi = Lp' = p$ is a weak equivalence as desired.

    \end{proof}

\begin{lemma}
Let $J'$ be the set of generating trivial cofibrations of $\sSetB$, $D\Sd^2J$. Then $i\colon  \sSetBC \rightarrow \sSetB$ takes $cell(LJ')$ to weak equivalences.
\label{secondcond}
\end{lemma}
\begin{proof}

We have previously shown $D\Sd^2\Delta[n]$ and $D\Sd^2\Lambda^k[n]$  are objects of $\sSetBC$ (Theorem \ref{thmbc}). Furthermore, the subdivided horn includes fully into the subdivided $n$-simplex with those vertices; this is since any simplex in $\Sd^2\Delta[n]$ is either also a simplex in one of its subdivided faces, or contains a vertex in its interior. Then, let $j'\colon X \to Z$ be a generating trivial cofibration in $\sSetB$; then the map to the pushout $Y \to Z \sqcup_X Y$ (when taken in $\sSetB$) is a weak equivalence.

 By Lemma $\ref{mainlemma}$, $i$ takes pushouts along $X \to Z$ to weak equivalences in $\sSetB$. Furthermore, $i$ preserves directed colimits as established in the proof of Lemma $\ref{slemma}$. 
 Hence, $i$ maps the transfinite compositions of pushouts of maps in $LJ'$ to transfinite compositions of homotopy pushouts of maps in $J'$ in $\sSetB$, ie. $i(cell(LJ')) \simeq cell(LJ')$. The maps in $cell(LJ')$ are  weak equivalences in $\sSetB$ since they are mapped by $\Ex^2i$ to weak equivalences in $\sSet$ \cite[Theorem 1.2]{fjellbonsSet}. Thus $i$ brings all maps in $cell(L J')$ to weak equivalences.

\end{proof}
  Next, we will show some results which will be needed to show the Quillen equivalence once we have established a Quillen adjunction.

  \begin{definition}
  \label{localdef}
        Let $\Sd^2\partial\Delta[n] \hookrightarrow \Sd^2\Delta[n]$ be a generating cofibration, $f \colon \Sd^2\partial\Delta[n] \to Y$ any map of simplicial sets such that $Y$ has Property $C$, and denote their pushout in $\sSetB$ as $P =D(\Sd^2\Delta[n] \sqcup_{ \Sd^2\partial\Delta[n]} Y)$, and in $\sSetBC$ as $LD(\Sd^2\Delta[n] \sqcup_{ \Sd^2\partial\Delta[n]} Y)$. 
        Fix an $n$-simplex $x \in \Sd\Delta[n]$, and let $X_x \subseteq P$  denote $\Sd(x) \subseteq \Sd^2\Delta[n]$. 

We refer to the intersection of the image in $P$ of a subcomplex of $\Sd^2\Delta[n]$ with $\Sd^2\partial\Delta[n]$ as the \textit{attaching boundary} of that subcomplex. And given $X_x$, we refer to the subcomplex generated by simplices with a face on the attaching boundary as $K_x$.

  \end{definition}
   \begin{lemma}
       If $Y$ has Property $C$, $P$ cannot contain two distinct nondegenerate simplices with the same vertices if all of those vertices lie in the attaching boundary. \label{diffboundary}
   \end{lemma}
   \begin{proof} 
   The map $Y \to P$ is full by Lemma \ref{lemma2}. Hence, every simplex of $P$ whose vertices lies in the attaching boundary belong to $Y$. Two such nondegenerate simplices with the same vertices must therefore be equal since $Y$ has Property C.

   \end{proof}
   
  We say that a simplicial subset $U\subseteq P$ is closed under $L$ if it is closed under the identifications defining $\eta_P:P\to iLP$; if $\sigma\in U$ and $\eta_P(\sigma)=\eta_P(\tau)$, then $\tau\in U$.
   \begin{lemma}
  The subcomplexes $X_x,K_x$, and $Y$ are closed under $L$.
   \label{saturationremark}
   \end{lemma}
   \begin{proof}

 Let the $q$-simplex $\sigma \in X_x$ be represented as $F_0\subsetneq\cdots\subsetneq F_q$, a chain of faces in $x$. Suppose that $\sigma$ is identified with another $q$-simplex $\tau$ by $L$. If $\sigma$ has an interior vertex, then $\tau$ has the same interior vertices and can be represented as a chain with maximal face $F_q$; every face in the chain is a face of $F_q$ and hence in $x$, so $\tau \in X_x$. If $\sigma$ has no interior vertex, then all its vertices lie in the attaching boundary, so Lemma $\ref{diffboundary}$ implies that $\sigma=\tau$ in $P$, and so $\tau \in X_x$.
 Since every simplex of $X_x$ identified by $L$ belongs to $K_x$, the subcomplex $K_x$ is also closed under identifications by $L$.
 
 Finally, $Y\to P$ is full and $Y$ has Property $C$, so $Y$ is closed as well.
   \end{proof}
   \begin{corollary}
\label{localinclusions}
The canonical maps
$LX_x\to LP,$ $LK_x\to LP,$ and $LY\to LP$ are simplicial inclusions with images
$\eta_P(X_x),$ $ \eta_P(K_x),$ and  $\eta_P(Y)$, respectively.
\end{corollary}

   Closure under $L$ allows subsimplicial sets to participate in particular pushouts in $\sSet$ or $\sSetB$ coinciding with those taken in $\sSetBC$.   
   \begin{lemma}
\label{closedunion}
Let $U,V\subseteq P$ in $\sSetB$ be simplicial subsets such that $P=U\cup V$. Suppose that $U$ and $V$ are closed under the equivalence relation defining $\eta_P\colon P\to LP$. Then $U\cap V$ is also closed, the canonical maps $LU\to LP,
LV\to LP,$ and $L(U\cap V)\to LP$ are simplicial inclusions, and the square
\[
\begin{tikzcd}
L(U\cap V) \arrow[r] \arrow[d] & LV \arrow[d]\\
LU \arrow[r] & LP
\end{tikzcd}
\]
is a pushout in $\sSet$

\end{lemma}

\begin{proof}
Since $U$ and $V$ are closed under $L$, so is $U \cap V$; thus the canonical maps $LU\to LP$, $LV\to LP$, and $L(U\cap V)\to LP$ are inclusions with images $\eta_P(U), \eta_P(V)$, and $\eta_P(U\cap V)$ respectively.

Since $P=U\cup V$ and $\eta_P$ is surjective, $LP=\eta_P(U)\cup\eta_P(V)$.
Closedness and the identification of $L(U), L(V)$ with their images under $\eta_P$ in $L(P)$ also gives
$L(U)\cap L(V)=L(U\cap V)$.
Therefore $LP\cong LU\sqcup_{L(U\cap V)}LV$ taken in $\sSet$, as required.
\end{proof}

   \begin{lemma}
   \label{localbehavior}
Using the notation of Lemma \ref{saturationremark}, and denoting $K =  \bigcup_{x}K_x$ the quotient map $K \to LK$ is a weak equivalence.
   \end{lemma}
   \begin{proof}
  This uses essentially the same argument as Lemmas \ref{defret} and $\ref{defret2}$. Instead of showing a deformation retraction of the quotiented simplex onto the quotiented horn, however, we will show $K$ deformation retracts onto the quotiented boundary. Denote $N_x \subseteq x$ as the subcomplex of $x$ generated by simplices with a face on the boundary, and $N = \bigcup_x N_x$ so that the image of $N$ in $P$ is $K$. Consider a small ball $B^n$ centered around the barycenter of the entire $|\Sd^2\Delta[n]|$, and take the radial retraction $H$ of $|\Sd^2\Delta[n]|\backslash B^n$ to $|\Sd^2\partial\Delta[n]|$. We may restrict this deformation retraction to $\tilde{H}$ on $N$: each ray intersects some $N_x$ in a single interval ending on the attaching boundary, and $C_x := \overline{|\Sd(x)|\backslash |N_x|}$ is convex. To see this, let $c$ denote the barycenter of $\Delta[n]$ and write the vertices of $x$ as $c,v_1,\ldots,v_n$. If
$p=\lambda_c c+\lambda_1v_1+\cdots+\lambda_nv_n$, since $C_x$ are the simplices of $\Sd(x)$ containing $c$, the barycentric subdivision gives
$|C_x|=\{p\in|x|:\lambda_c\geq\lambda_i\text{ for every }i\}$. It is an intersection of closed affine half-spaces and therefore convex. 
Each ray starting at $c$ thus meets $C_x$ at an initial segment and exits at a unique point, and crosses into $N_x$ via a terminal segment intersecting the attaching boundary. Hence the radial retraction restricts to each $N_x$ as well as on $N$. 

  Denote the image of $\Sd^2\partial\Delta[n]$ as $F \subseteq P$. Since only simplices on the boundary may be identified by $D$, we may use Lemma \ref{defret} to show this deformation retraction $\tilde{H}\colon |N| \times I \to |N|$ descends to a deformation retraction $|K| \times I \to |K|$ to the quotiented boundary. Hence $|F| \to |K|$ is a weak equivalence. Futhermore, we may use the argument of Lemma \ref{defret2}, by similarly using $\tilde{H}$ and accounting for quotienting in simplices with more than one vertex on the boundary (since these are a subcomplex of $K$), to define a deformation retraction $|LK| \times I \to |LK|$ to the quotiented boundary $LF$ (the image of $\Sd^2\partial\Delta[n]$ in $LP$). Hence we have that $|LF|\to |LK|$ is a homotopy equivalence. The maps are hence weak equivalences in $\sSetB$.
  By naturality of $\eta$ we thus have the commuting diagram:
  \[
  \begin{tikzcd}
      F \arrow[r, "\sim"]\arrow[d, "\eta_F"] & K\arrow[d, "\eta_K"]\\\
      LF \arrow[r, "\sim"] & LK
  \end{tikzcd}
  \]
  
Finally, $\eta_F\colon F \to LF$ is a weak equivalence as a direct consequence of Lemma \ref{diffboundary}; indeed, $F$ may be considered a subset of $Y$ and hence has Property C. Thus $\eta_K \colon K \to LK$ is a weak equivalence.

   \end{proof}

   We want to glue together compatible pieces to form the pushout, which will use the following gluing lemma for proper model categories \cite[pg 7]{fjellbonsSet}, \cite[pg 246-247]{HirschhornModelCategories}.
   \begin{lemma}[Gluing Lemma (Hirschhorn)]
       In a proper model category, consider commutative diagrams of the following form:
       \[
\begin{tikzcd}
    A \arrow[d, "\sim"] &\arrow[l] B \arrow[d, "\sim"]\arrow[r]&C\arrow[d, "\sim"] \\
        A'  &\arrow[l] B' \arrow[r]&C'
\end{tikzcd}
\]
If all three vertical maps are weak equivalences and at least one map in each row is a cofibration, the induced map on pushouts $A \sqcup_B C \to A' \sqcup_{B'} C'$ is also a weak equivalence.
   \end{lemma}
\begin{lemma}
    \label{gluingtogether}
    Let $Y\in\sSetB$ have Property C,
let $f\colon \Sd^2\partial\Delta[n]\to Y$,
and let $P=D(\Sd^2\Delta[n]
\sqcup_{\Sd^2\partial\Delta[n]}Y)$.
Then the unit $\eta_P\colon P\to iLP$ is a weak equivalence.
\end{lemma}
\begin{proof}
Denote $C$ as the contractible subcomplex of interior simplices of $P$ which do not contain any faces on the attaching boundary. Note that equivalently, this is the subcomplex of $n$-simplices which contain the vertex common to all local subcomplexes.
As before, let $K = \bigcup_xK_x$. We would like to use the gluing lemma on $P = C \sqcup_{K \cap C} (Y\cup K)$. We know $K \cap C$ is a subcomplex which is fully contained in the interior and hence not identified by $L$ or $D$; therefore $K \cap C \to L(K \cap C)$ is a weak equivalence, and similarly $C \to LC$ is also. By Lemma \ref{localbehavior}, we have $K \to L(K)$ is a weak equivalence. Since $L$ does not identify any simplices contained entirely in $Y$ or the attaching boundary ($Y \cap K$), $Y \to LY$ and $(Y \cap K) \to L(Y \cap K)$ are weak equivalences. Furthermore, $K, Y, Y\cap K$ are closed under identifications by $L$, since each of $K_x$ and $Y$ is closed under $L$, hence their unions and intersections are. We then may apply the gluing lemma to obtain that $(Y \cup K) \to L(Y \cup K)$ is a weak equivalence.  Then $Y \cup K$, is also closed under identifications by $L$, as well as $C$ and $K \cap C$, so Lemma \ref{closedunion} implies that $LP = LC \sqcup_{L(K \cap C)} L(Y \cup K)$.
 \[
\begin{tikzcd}
   (Y \cup K) \arrow[d, "\sim"] &\arrow[l] K \cap C \arrow[d, "\sim"]\arrow[r]&C\arrow[d, "\sim"] \\
        L(Y \cup K)  &\arrow[l] L(K \cap C) \arrow[r]& LC
\end{tikzcd}
\]
After realization, the attaching maps are inclusions of CW subcomplexes and hence cofibrations. The gluing lemma therefore shows that $P\to iLP$ is a weak equivalence.
\end{proof}
        
    \begin{lemma}
    For every cofibrant object $Z \in \sSetB$, the unit
    $$\eta_Z\colon Z \to iLZ$$ is a weak equivalence.
    \label{transfiniteinduction}
    \end{lemma}
    \begin{proof}
        Let $D\Sd^2I\colon \{D\Sd^2\partial\Delta[n] \to D\Sd^2\Delta[n]|n\geq 0\}$ be the generating cofibrations of $\sSetB$. Suppose that $X \in cell(D\Sd^2I)$, then we may write $X$ as a transfinite sequence $$X_0\to X_1\to\cdots\to X_\alpha\to\cdots$$ in which $X_0=\emptyset$, $X_i \to X_{i+1}$ is obtained by a pushout along a generating cofibration, and at every limit ordinal $\lambda$, $X_\lambda=\colim_{\alpha<\lambda}X_\alpha$. We prove by transfinite induction that $\eta_{X_\alpha}$ is a weak equivalence. The result is immediate for $X_0=\emptyset$.

        Let $A\to B$ be a map in $D\Sd^2I$, and let $X_{\alpha+1} =X_\alpha\sqcup_A B$. Suppose that $\eta_{X_\alpha}\colon X_{\alpha} \to iLX_\alpha$ is a weak equivalence. Now consider the following composition of pushout squares:
        \[
        \begin{tikzcd}
            A\arrow[r] \arrow[d]&X_\alpha\arrow[d] \arrow[r, "\eta_{X_\alpha}"] & iLX_\alpha \arrow[d]\\
            B\arrow[r]& X_{\alpha+1} \arrow[r, "q_\alpha"] & Q_{\alpha}
        \end{tikzcd}
        \]
         By composing the pushout $q_\alpha$ obtained from $X_{\alpha+1} \sqcup_{X_\alpha} iLX_{\alpha}$ with the unit $\eta_{Q_{\alpha}}$, we obtain the composite $X_{\alpha+1} \xrightarrow{q_\alpha} Q_\alpha \xrightarrow{\eta_{Q_\alpha}} iLQ_{\alpha}$. Since $X_\alpha \to X_{\alpha+1}$ is a cofibration and $\eta_{X_\alpha}$ by the inductive hypothesis $\eta_{X_\alpha}$ is a weak equivalence, we may use left properness of $\sSetB$ to obtain that $q_\alpha$ is also a weak equivalence. Lemma \ref{gluingtogether} gives that unit $\eta_{Q_\alpha}\colon Q_\alpha \to iLQ_\alpha$ is a weak equivalence. To show $iLq_\alpha$ is an isomorphism and hence a weak equivalence, it suffices to show $Lq_\alpha$ is. Since $L$ is left adjoint, $Lq_\alpha$ is the map induced on the pushout by:
         \[
        \begin{tikzcd}
            LB\arrow[r, leftarrow] \arrow[d, "id"]&LA\arrow[d, "id"] \arrow[r] & LX_\alpha \arrow[d, "L\eta_{X_\alpha}"]\\
            LB\arrow[r, leftarrow]& LA \arrow[r] & LiLX_{\alpha}
        \end{tikzcd}
        \]
         The right map is an isomorphism by the triangle identities: $id_{LX_\alpha} = \varepsilon_{LX_\alpha}\circ L\eta_{X_\alpha}$, and $\varepsilon_{LX_\alpha}$ is an isomorphism since $L\dashv i$ defines a reflective subcategory, so $L\eta_{X_\alpha}$ is also. Hence so is their induced map on the pushouts. Then by 2-of-3, $\eta_{X_{\alpha+1}} \colon X_{\alpha+1} \to iLX_{\alpha+1}$ is a weak equivalence.

Now let $\lambda$ be a limit ordinal. Since $L$ and $i$ (by the proof of Lemma \ref{slemma}) preserve directed colimits, $\eta_{X_\lambda} \cong \colim_{\alpha<\lambda}\eta_{X_\alpha}$. Weak equivalences in $\sSetB$ are preserved under directed colimits: they are detected by $Ex^2i$, which preserves directed colimits, and weak equivalences of simplicial sets are preserved under directed colimits. Thus $\eta_{X_\lambda}$ is a weak equivalence.

Finally, any cofibrant object $Z$ can be written as a retract of $X \in cell(D\Sd^2I)$, so by naturality of the unit, $\eta_Z$ is also a retract of $\eta_X$. Then $\eta_Z$ is a weak equivalence since retracts are closed under weak equivalences.
     \end{proof}
   
We are now ready to show:
\begin{theorem}
\label{sSetBCmodel}
The  model structure on $\sSetB$ of Theorem \ref{thm122} can be right-transferred to $\sSetBC$ along the adjunction  \[
        \begin{tikzcd}
            \sSetB \arrow[r, shift left=.75ex,"L"{name=G}] & \sSetBC\arrow[l, shift left=.75ex, "i"{name=F}] 
            \arrow[phantom, from=F, to=G, "\scriptstyle\dashv" rotate=270]      
        \end{tikzcd}
    \] Moreover, this Quillen pair is a Quillen equivalence.

\end{theorem}

\begin{proof}
This is an application of Theorem \ref{thm1}. Since $\sSetBC$ is a reflective subcategory of $\sSetB$, which is closed under all small limits and colimits, $\sSetBC$ is also bicomplete. The first condition follows directly from Lemma \ref{slemma}, since the domains of $LI$ and $LJ$ are subdivided horns or boundaries, which are finite and hence small; and Lemma \ref{secondcond} is the second condition. This proves the existence of the right-transferred model structure, and as a consequence, we have that $(L, i)$ is a Quillen pair. To show it is a Quillen equivalence,  let $f$ be a map of simplicial sets in $\sSetBC$. Since $f$ is a weak equivalence if and only if $i(f)$ is a weak equivalence, the Quillen pair $(L, i)$ is a Quillen equivalence if and only if for all cofibrant $X \in \sSetB$, the unit $X \to iL X$ is a weak equivalence in $\sSetB$ \cite[Lemma 3.3]{Erdal_2019}, which is given by Lemma \ref{transfiniteinduction}. 
\end{proof}

\begin{corollary}
    There is a chain of Quillen equivalences:
    \[
    \begin{tikzcd}
            \sSet \arrow[r, shift left=.75ex,"D\Sd^2"{name=G}] & \sSetB\arrow[l, shift left=.75ex, "\Ex^2i"{name=F}] 
            \arrow[phantom, from=F, to=G, "\scriptstyle\dashv" rotate=270]\arrow[r, shift left=.75ex,"L"{name=H}] & \sSetBC\arrow[l, shift left=.75ex, "i"{name=I}]\arrow[phantom, from=H, to=I, "\scriptstyle\dashv" rotate=270]
            \end{tikzcd}\]
            In particular, the model structure on $\sSetBC$ of Theorem $\ref{sSetBCmodel}$ is Quillen equivalent to the standard model structure on simplicial sets.
\end{corollary}

\section{From Simplicial Sets to Ordered Simplicial Complexes}
The objects of the category $\sSetBC$ essentially have the properties which distinguish simplicial complexes, but they are still objects of $\sSet$. We want to establish a model structure on a category whose objects are defined directly from the definition of simplicial complexes, independent of simplicial sets.

In this section we define $\Osc$, the category of ordered simplicial complexes. We will also construct the left adjoint to the canonical functor \cite[Proposition 10.2.2]{finitespaces} from $\Osc$ to $\sSetBC$.

\begin{definition}
An \textit{ordered simplicial complex} $(S$, $K \subseteq \mathcal{P}(S), \mathcal{R})$ is a set of vertices $S$, a collection of subsets $K$, and a relation $\mathcal{R}$ such that: (1) $(S, K)$ is a simplicial complex, (2) $\mathcal{R}$ is transitive, antisymmetric, and reflexive on vertices (ie. a partial order), and (3) $\forall T \in K$, $\mathcal{R}|_T$ is a total order: $\forall t_1, t_2 \in T$, either $t_1 \leq t_2$ or $t_2 \leq t_1$, and if $t_1 \leq t_2$ and $t_2 \leq t_1$ then $t_1 = t_2$.
\end{definition}

Let $\Osc$ be the category whose objects are ordered simplicial complexes and morphisms are order-preserving maps defined on vertices.

  \begin{definition}
 The functor $U\colon \Osc \rightarrow \sSetBC$ sends an ordered simplicial complex $Y$ to the simplicial set $U(Y)$, with $U(Y)_0$ the vertices of $Y$. We get a nondegenerate $n$-simplex in $U(Y)$ for each subset of size $n$ in $Y$; we freely add the degenerate simplices as well as the degeneracy maps by letting $s_i(v_0,\dots,v_n) = (v_0,\dots,v_i,v_i,\dots,v_n)$.
 
 Furthermore, for each simplex, define the face maps as $d_i(v_0,\dots,v_n) = (v_0,\dots,v_{i-1},v_{i+1},\dots,v_n)$, which is an $(n-1)$-simplex since $Y$ is downward closed by inclusion.
 
 A map between ordered simplicial complexes $g\colon  Y \rightarrow Y'$ is sent to $Ug$ which is the simplicial map sending each nondegenerate $n$-simplex $(a_1, \dots, a_n)$ in $U(Y)$, to $(g(a_1), \dots, g(a_n))$ in $U(g(Y))$. This fully determines the simplicial map, since any map $f\colon  X \rightarrow Y$ in $\Osc$ is determined completely by the map between vertices, and so is the map $Uf\colon  U(X) \rightarrow U(Y)$ as a morphism in $\sSetBC$ by Proposition \ref{vertexmap}.
 \end{definition}

 Note that the simplices of an ordered simplicial complex are totally ordered, which allows one to define face and degeneracy maps without making an arbitrary choice. This would not be possible with unordered simplicial complexes.
 
 Although ordered simplicial complexes and simplicial sets with Properties B and C are closely connected, the relation between the two categories is much more subtle than initially believed in \cite[Theorem 12.1.8]{finitespaces}, as we detail in the following remarks.
\begin{remark}
\label{rmkloop}
The category $\sSetBC$ is not equivalent to $\Osc$. This seems to be in conflict with statements made in \cite{finitespaces}. Note that in \cite[Theorem 12.1.8]{finitespaces}, the following is claimed: ``A simplicial set $X$ is a simplicial complex, (ie. the simplicial realization of an ordered simplicial complex), if and only if it satisfies Properties B and C". But this reformulation of a simplicial complex does not exactly hold since it seems only the converse is true. It can be easily checked that any simplicial complex in $\Osc$ has a simplicial set representation with Properties B and C after applying the canonical injection functor $U$. However, there may be simplicial sets with Properties B and C which do not arise directly from ordered simplicial complexes, such as the simplicial set consisting of the three nondegenerate 1-simplices depicted below:
\[
\begin{tikzcd}  & b \arrow[dr] \\a \arrow[ur] & & c \arrow[ll]
\end{tikzcd}
\]
This could not be the image of an ordered simplicial complex because a partial order which restricts to these orderings for the simplices would cause all the vertices to be identified.

     To elaborate: $U$ has been identified as the canonical functor from $\Osc$ to $\sSet$, and hence to $\sSetBC$, by sending an ordered simplicial complex to the obvious choice of simplicial set with the same nondegenerate simplices \cite{finitespaces}. However, $U$ is not in fact fully faithful as previously assumed in Proposition 10.2.2 of \cite{finitespaces}. Let $X \in \Osc$ consist of two vertices $a$ and $b$, no other simplices, and the partial order relation $a \leq b$. Then there are 3 elements in $\Osc(X,X)$ but 4 elements in $\sSetBC(UX,UX)$; in particular, there is no map swapping $a$ and $b$ in $\Osc$, but there is in $\sSetBC$.

     Thus, the relationship between $\Osc$ and $\sSetBC$ is more subtle, and $U$ is not an embedding of $\Osc$ into $\sSet$ or $\sSetBC$. Rather, the essential image of $U$ consists of simplicial sets such that the order relation generated by taking the natural order on each simplex is still a partial order. Hence we can consider two 0-simplices in a simplicial set to be comparable if and only if they are connected by a chain of 1-simplices. Given any simplicial set $X$ in $\sSetBC$, we would like to find a representative $X/\sim$ in the essential image, which we obtain by identifying loops of 1-simplices.
\end{remark}

\begin{definition}\label{defn:quotient0}
    Let $X \in \sSetBC$, and define $p_0: X_0 \rightarrow X_0/\sim_0$ to be the projection of vertices to the quotient defined by contracting any loops existing in the 1-skeleta of $X$. Equivalently, a vertex $x$ is identified with $y$ if and only if there exists an opposing pair of morphisms $x\to y$ and $y\to x$ in the category $\pi_1 X$. 
\end{definition}

\begin{definition}
\label{simdef}
    Given $X \in \sSetBC$, define $X/\sim$ as the simplicial set obtained by the pushout   \[    \begin{tikzcd}
X_0 \arrow[hookrightarrow]{r}{i} \arrow[swap]{d}{p_0} & X \arrow{d}{p} \\%
X_0/\sim_0 \arrow{r}{j}& X/\sim
\end{tikzcd} \]
taken in $\sSetBC$, where $i$ is the inclusion into $X$ and $p_0$ is the projection of the discrete simplicial set with vertices $X_0$ to the discrete simplicial set with vertices $X_0/\sim_0$.

Let the induced map $p\colon  X \rightarrow X/\sim$ be denoted as the projection to the quotient of $X$.
\end{definition}
\begin{proposition}
 $X/\sim$ as defined in Definition \ref{simdef} can be explicitly constructed as the following simplicial set:
\begin{itemize}
    \item The 0-simplices $(X/\sim)_0 = X_0/\sim_0$ as in Definition \ref{defn:quotient0}.
    \item For each $n\geq 1$, given an ordered list of vertices $([a_0],\dots,[a_n])$, $X/\sim$ has a unique $n$-simplex with those vertices if and only if there is an $n$-simplex $(b_0,\dots,b_n)$ in $X$ such that $[b_i]=[a_i]$ for all $i$.
    \item The face maps are given by \[d_i([a_0],\dots,[a_n]) = ([a_0],\dots,[a_{i-1}], [a_{i+1}], \dots, [a_n]).\] Note that if $(a_0,\dots,a_n)\in X_n$ is a simplex witnessing the existence of $([a_0],\dots,[a_n])\in (X/\sim)_n$, then $d_i(a_0,\dots,a_n)$ will witness the existence of $([a_0],\dots,[a_{i-1}], [a_{i+1}], \dots, [a_n])\in (X/\sim)_{n-1}$; hence the face maps are well-defined.
    \item The degeneracy maps are given by $s_i([a_0],\dots,[a_n]) = ([a_0], \dots,[a_{i}], [a_{i}], \dots,[a_{n+1}])$; these are well-defined by an analogous argument. 
\end{itemize}

The projection map $p\colon  X \rightarrow X/\sim$ maps a vertex $a\in X_0$ to $p_0(a)=[a]$, and an $n$-simplex $(a_0,\dots,a_n)\in X_n$ to the unique $n$-simplex $([a_0],\dots,[a_n])\in (X/\sim)_n$.

\end{proposition}
\begin{proof}
    Suppose we are given an arbitrary $Y \in \sSetBC$ and maps $f\colon  X \rightarrow Y$ and $g\colon  X_0/\sim_0 \rightarrow Y$ such that $fi = gp_0$. Since $fi(v) = gp_0(v)$ for all $v \in X_0$, and $i$ is the identity on $X_0$ (hence surjecting onto the vertices of $X$), $f$ on vertices can be determined from $gp_0(v)$. However, a map in $\sSetBC$ is fully determined by the map on vertices by Proposition \ref{vertexmap}; hence $f$ is uniquely determined by $g$. The induced map $q\colon  (X/\sim) \rightarrow Y$ maps $[v] \in (X/\sim)_0$ to $fi(v) = f(v)$, which commutes with the diagram since $p(v) = [v]$; as a map in $\sSetBC$, this fully determines the unique map $q$.
\end{proof}

Now we use the previous construction to define a left adjoint to $U$.
 \begin{definition}
 \label{defquot}
 The functor $F\colon  \sSetBC \rightarrow \Osc$ sends a simplicial set $X \in \sSetBC$ to the ordered simplicial complex with vertices the set  $X_0/\sim_0$, and a (simplex) subset $K$ for each nondegenerate simplex of $X$ whose vertices all have distinct image in $X_0/\sim_0$. Each subset $T \in K$ is given the canonical order relation of the corresponding simplex in $X$, and these are extended to the order relation $\mathcal{R}$ on the vertices by taking the transitive closure of the union of all these relations. 

 A map between simplicial sets $f\colon  X \rightarrow Y$ is sent to a map $Ff\colon  FX \rightarrow FY$ which takes a vertex $[a] \in (X/\sim)_0$ to $[f(a)]$.
 \end{definition}
 \begin{example}
We give an example of the construction of $FX$ for a simplicial set $X \in \sSetBC$. Suppose that $X$ is the simplicial set in $\sSetBC$ consisting of 6 0-simplices, 7 nondegenerate 1-simplices, and 1 nondegenerate 2-simplex as depicted in the following image:
    \[\includegraphics[width=0.3\linewidth]{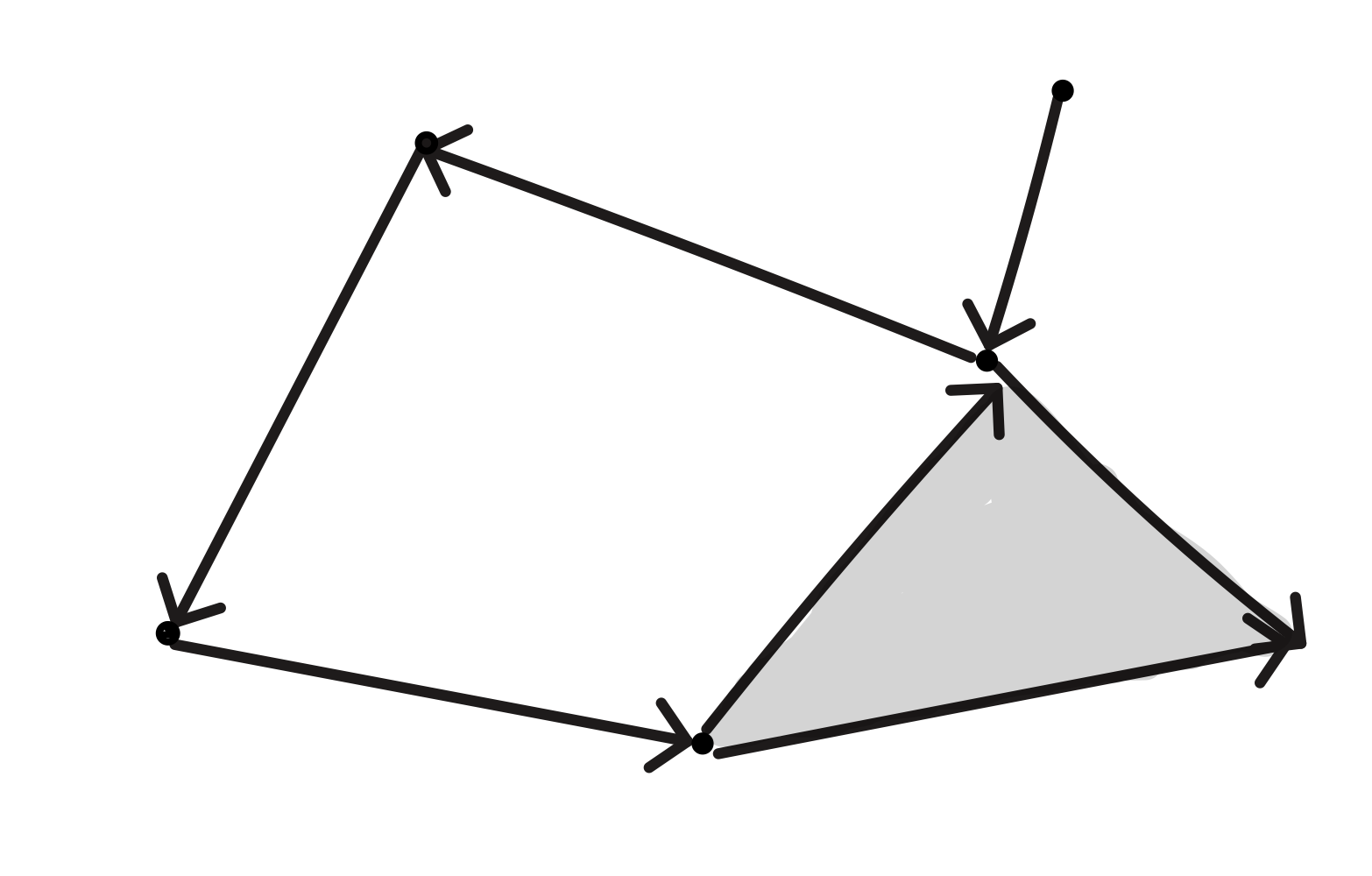}\] 
    The quotient $X/\sim$ is given by collapsing the left loop of 4 1-simplices to a point, which would cause a face of the 2-simplex, and hence the 2-simplex itself to become degenerate (since we are in $\sSetBC$). The result is that $X/\sim$ is the simplicial set given by 3 0-simplices $\{a, b, c\}$ and nondegenerate 1-simplices $a \to b$ and $ b \to c$. Then the ordered simplicial complex $FX$ obtained by taking these vertices, simplices, and order relations generated by the order relations on the simplices of $X/\sim$ is the vertex set $\{a, b, c\}$, simplices $\{(a), (b), (c), (a,b), (b,c)\}$, and order relation $a < b <c$.
 \end{example}

  \begin{remark}
 The image of $F$ is antisymmetric, and hence a valid partial order: each simplex of $FX$ corresponds to a unique simplex of $X$ with distinct vertices in $(X/\sim)_0$. Equivalently, each simplex of $FX$ corresponds to the unique simplex of $(X/\sim)$ with the same distinct vertex representatives. Since there are no loops in the 1-skeleta of $(X/\sim)$, there are no opposing relations formed by passing to the transitive closure; hence, the resulting relation is antisymmetric and thus a partial order. 

 Furthermore, we would like to think of $X/\sim$ as the ``largest" simplicial set in $\sSetBC$ obtained by identifying vertices of $X$ in which the transitive closure of the order relations on simplices is a partial order. Indeed, any other quotient of $X$ obtained by identifying vertices, whose transitive closure is a partial order, factors through $X \rightarrow X/\sim$.
 
 \end{remark} 
\begin{proposition}
We have an adjunction:
\[
        \begin{tikzcd}
            \sSetBC \arrow[r, shift left=.75ex,"F"{name=G}] & \Osc\arrow[l, shift left=.75ex, "U"{name=F}] 
            \arrow[phantom, from=F, to=G, "\scriptstyle\dashv" rotate=270]     
        \end{tikzcd}
    \]
    \label{oscad}
\end{proposition}
\begin{proof}
    Given an ordered simplicial complex $Y$, note that $FUY$ is the ordered simplicial complex with the same simplices as $Y$, but there are no order relations in $FUY$ which are not in $Y$. Then, let the counit $\varepsilon\colon  FU \to 1_{\Osc}$ be the inclusion.

    For the unit $\eta: 1_{\sSetBC} \rightarrow UF$, first note that $UFX \cong X/\sim$; we then define $\eta_X$ as the projection to the quotient $X \to X/\sim$. 

Now we need to show the following diagrams are natural: \\
    \[    \begin{tikzcd}
X \arrow{r}{\eta_X} \arrow[swap]{d}{f} & UFX \arrow{d}{UFf} \\%
X' \arrow{r}{\eta_{X'}}& UFX'
\end{tikzcd} \quad \quad \begin{tikzcd}
FUY \arrow{r}{\varepsilon_Y} \arrow[swap]{d}{FUg} & Y \arrow{d}{g} \\%
FUY' \arrow{r}{\varepsilon_{Y'}}& Y'
\end{tikzcd}
\]
$\\$
Let $X$ be a simplicial set in $\sSetBC$ and $Y$ an ordered simplicial complex in $\Osc$. Since $\eta_X$ is the induced map $p\colon  X \rightarrow X/\sim$, naturality of the first square follows from naturality of the pushout. For the second square, it suffices to note that $\varepsilon_Y$ is the identity on vertices, and that $FUY$ and $Y$ have the same vertices and simplices.

 Now we verify the triangle identities:
\[\begin{tikzcd}[row sep=3em]
FX \arrow{r}{F\eta_X} \arrow[""{name=foo}]{dr}[swap]{1_{FX}} & FUFX \arrow{d}{\varepsilon FX} \arrow[Rightarrow, from=foo, swap, near start, "="]& UY \arrow{d}[swap]{\eta UY} \arrow[""{name=bar, below}]{dr}{1_{UY}} & \\
& FX & UFUY \arrow{r}[swap]{U\varepsilon_Y} \arrow[Rightarrow, to=bar, swap, near start, "="] & UY
\end{tikzcd}\]
An ordered simplicial complex $FX$ is sent to a simplicial set with the same simplices when we apply $U$, and because all the order relations on vertices must be compatible with the ordering of the $FX \in \Osc$, applying $F$ gives the same ordered simplicial complex. Hence $F\eta_X$ is the identity. The inclusion of the identical ordered simplicial complex is also the identity, hence $\varepsilon FX\circ F\eta_X = 1_{FX}$. Given an $Y\in \Osc$, $UY$ is a simplicial set with the simplices of $Y$ and order relations compatible with $Y$. Then applying the unit $\eta$ projects $UY$ onto itself, so $\eta UY$ is the identity. The inclusion of $FUY$ into $Y$ is the inclusion of an ordered simplicial complex with the same simplices and no extra order relations; whiskering with $U$, gives identical simplicial sets, so $U\varepsilon_Y$ is the identity on the resultant simplicial sets. Hence $U\varepsilon_Y\circ\eta UY = 1_{UY}$.
\end{proof}

\section{Equivalence of Homotopy Theories for Ordered Simplicial Complexes and Simplicial Sets}
As mentioned in the previous section $\Osc$ is not a reflective subcategory of $\sSetBC$, the canonical functor $U\colon  \Osc \to \sSetBC$ is not fully faithful. However, the adjunction $(F, U)$ of Theorem \ref{oscad} enjoys special properties which allow us to lift the model structure from $\sSetBC$ of Theorem \ref{sSetBCmodel} to $\Osc$.
We will again use the transfer theorem, Theorem \ref{thm1}, to lift the model structure from $\sSetBC$ to $\Osc$. To start, we construct limits and colimits in $\Osc$ with the following lemma.
\begin{lemma}
 $\Osc$ has all small limits and colimits.
 \label{osccolimits}
 \end{lemma}
 \begin{proof}
For limits, it suffices to show that $\Osc$ has all products and equalizers. For products $\prod J_i$ of $J_i = (S_i,  K_i,\mathcal{R}_i) \in \Osc$, define the vertices as $$S_{\prod J_i}:= \prod S_i\,,$$ which is a set of tuples of vertices, where each tuple is composed of one vertex for each $J_i$. This set has a partial order in which each index inherits the partial order from the respective ordered simplicial complex, and we have a relation if and only if all indices are comparable. Now define the set of simplices $$K_{\prod J_i} := \{\sigma = \text{sets of totally ordered subsets of } S_{\prod J_i}|p_i(\sigma) \in K_i\}$$
Note that simplices are not constrained by set size in $\Osc$, so we can take arbitrarily large products. This is a generalization of the construction in \cite{bergner2024simplicialsetstopologycategory} which defines the product between two ordered simplicial complexes.

For equalizers,  given a diagram $J\colon  (A \rightrightarrows B)$ in $ \Osc$, we first find the equalizer of $UJ$ in $\sSetBC$. This exists (as $\sSetBC $ is complete, as established in Theorem \ref{sSetBCmodel}) and agrees with the equalizer in $\sSet$, as $\iota\colon\sSetBC\to\sSet$ is a right adjoint; let's denote this $X\hookrightarrow UA$.

We define the ordered simplicial complex $C$ by giving it the same vertices as $X$. It remains to define an order on its vertices: we let $a \leq b \in C$ if $i(a) \leq i(b) \in A$. This can easily be checked to be the limit of $J$, as $C = \lim J$ inherits all order relations from $A$. In particular, $C$ has the maximal order relations such that $i$ is still a map. 
 
 To show that small colimits exist, it suffices to show that coproducts and coequalizers exist in $\Osc$. For coproducts, we need to show that given a collection $\{A_i\}_{i\in\mathcal{I}}$ in $\Osc$ for $\mathcal{I}$ an indexing set, there exists $B \in $ $\Osc$ such that Hom$_{\Osc}(B, X) \cong \prod_{i\in\mathcal{I}}$ Hom$_{\Osc}(A_i,X)$. We let $B$ be the disjoint union of all the $A_i$; the isomorphism is given by the map sending each tuple $(f_i)_{i \in \mathcal{I}}$ on the right to the corresponding map on the left with each disjoint $A_i$ as its domain. Furthermore, since the $A_i$ are disjoint, there are no further restrictions between components of the map, so that all possible combinations of the tuple on the right can be obtained by an appropriately defined map on the left; the reverse direction and injectivity follows immediately. 
 
 We define the following construction for coequalizers. Given a diagram $J=f, g:A \rightrightarrows B$ in $\Osc$, first consider the simplicial set $Z \in \sSetBC$ obtained by adding an edge $b_1 \to b_2$ to $UB$ if there is a relation $b_1 \leq b_2$ in $B$, and note that $UB$ includes naturally into $Z$. Now take the pushout $P$ of the diagram in $\sSetBC$:
 \[
\begin{tikzcd}
    \colim(UJ)_0 \arrow[d]\arrow[r]&\colim(UJ)\arrow[d] \\
    \colim(UA \rightrightarrows Z)_0\arrow[r, "h_0"]& P
\end{tikzcd} \]   $\colim(UA \rightrightarrows Z)_0$ cannot have fewer identifications of vertices, and so we can send a vertex in $\colim(UJ)_0$ to its equivalence class in $\colim(UA \rightrightarrows Z)_0$.  Then simplices in $FP$ do not conflict with the order relation on $B$, so for each $b_1 \leq b_2$ in $B$, we may add the order relations to the corresponding vertices of $FP$. This gives an ordered simplicial complex $X$ which we claim is $\colim(J)$, along with the map $h: B \to X$ defined via the map of vertices given by $h_0$. To see this, consider the coequalizing diagram, where $j\colon  B \to C$ is any morphism in $\Osc$ such that $jf = jg$.

The map $q$ can be constructed as follows: given $x \in X$, take $q(x) := j[h^{-1}(x)]$ where  $[h^{-1}(x)]$ is an element of the preimage of $x$ under $h$. In particular, this is well-defined; consider elements $b, b' \in [h^{-1}(x)]$, and $h(b), h(b')$ may be identified in $X$ in two ways, since $h$ was defined as a quotient of $\colim(UJ)$ from the order relations on the vertices from $B$. The first case is that there exists $a \in A$ such that $f(a)=b$ and $g(a) = b'$, and hence $jf(a) =j(b) = j(b') = jg(a)$. The second is from the order relations: that there are $a, a' \in A$ such that $f(a)\leq b \leq f(a')$, and $f(a) \leq b' \leq f(a')$, but $g(a') \leq g(a)$.
Since $j$ is an order-preserving map in $\Osc$, the condition $jf(a) = jg(a)$ in $C$ again requires $b, b'$ to be identified in the image of $j$. Furthermore, $q$ is unique, since a map of ordered simplicial complexes is determined fully by its action on vertices, and by construction the map $h$ is surjective. 

\end{proof}

 Similarly to Lemma \ref{slemma}, we can describe small objects in $\Osc$.

\begin{lemma}
An object in $\Osc$ is small if and only if it has a finite number of vertices, or equivalently, a finite number of simplices.
    \label{slemma2}
\end{lemma}
\begin{proof}
    The forward direction is the same argument for compactness of finite sets. For the reverse, suppose $X$ has a finite number of vertices, and we need to show $$\colim_j\Hom_\Osc(X, A_j)\cong \Hom_\Osc(X, \colim_j(A_j))$$ whenever $A$ is a filtered diagram in $\Osc$. Colimits in $\Osc$ were computed in Lemma \ref{osccolimits} by taking the colimit in $\sSetBC$, inheriting some quotient of its simplices, and adding in order relations between vertices. In addition, $\colim_j(UA_j)$ does not contain any loops (chains of directed 1-simplices starting and ending at the same vertex: see Definition \ref{loopdef}) since $n$-loops are small in $\sSetBC$ by Lemma \ref{slemma}, and since $F$ identifies vertices in a simplicial set in $\sSetBC$ only if it contains loops. Then since maps in $\Osc$ depend only on order-preserving maps between vertices, the condition to verify reduces to $$\colim_j\Hom_\OrdSet(X_0, (A_j)_0)\cong \Hom_\OrdSet(X_0, \colim_j((A_j)_0))$$ where $\OrdSet$ is the category of ordered sets and order-preserving maps. Given a map $f\colon  X_0 \to \colim_j((A_j)_0)$, $X_0$ is small in $\Set$ so we can factor $f$ as the unordered composite $f = gh$ for $h\colon  X_0 \to (A_k)_0$. However, since $f$ is order-preserving, for any order relation $x \leq x'$ in $X_0$, we have $f(x) \leq f(x')$ in $\colim_j((A_j)_0)$. Furthermore, there exists some $(A_l)_0$ mapping some $y \leq y'$ in $(A_l)_0$ to $f(x) \leq f(x')$ in $\colim_j((A_j)_0)$. 
    Then there exists some $l \geq k$ such that the composition of $h$ with $(A_k)_0 \to (A_l)_0$ is an order-preserving map commuting with all the other maps, and hence gives the required factorization.
\end{proof}

 \begin{remark}
   As in $\ref{thmbc}$, we will use the fact that $L D\Sd^2\Lambda^k[n] \cong \Sd^2\Lambda^k[n]$, and $ D\Sd^2\Delta[n] \cong \Sd^2\Delta[n]$, and refer to the generating trivial cofibrations, $\{D\Sd^2\Lambda^k[n] \hookrightarrow D\Sd^2\Delta[n]\mid n \geq 0, 0 \leq k\leq n$\}, as $\{\Sd^2\Lambda^k[n] \hookrightarrow \Sd^2\Delta[n]\mid n \geq 0, 0 \leq k\leq n$\} .
\end{remark}
\begin{lemma}
 \label{oldpushout} 
     Denote the pushout  in $\sSetBC$ along a generating (trivial) cofibration as $P =\Sd^2\Delta[n]\sqcup_{\Sd^2\partial\Delta[n]} W$, and $p$ the induced map $W \to P$. Then $P$ does not contain any loops if $W$ does not.
 \end{lemma}
 \begin{proof}
We show this for a generating cofibration, and the same follows for trivial cofibrations. Since simplices in subdivision are given by face inclusions, a directed chain of 1-simplices must increase in dimension going towards the interior of $\Sd^2\Delta[n]$. Since identifications of interior vertices cannot be made with vertices on the boundary or horn in the pushout, any loops formed must therefore occur completely on the image of the boundary $\Sd^2\partial\Delta[n]$ in the pushout. Since $p$ is an inclusion, we may identify $p(W)$ with $W$. But the image of the boundary $pf(\Sd^2\partial\Delta[n]) \subseteq p(W) = W$ must contain  a loop, so $W$ must also, which is false by assumption.
 \end{proof}
\begin{lemma}
 \label{eq2}
     Let $J$ denote the generating trivial cofibrations in $\sSetBC$. Then, $U$ takes the maps of $cell(FJ)$ to weak equivalences.
 \end{lemma}
 \begin{proof}
 First note $j \in J$ is a full simplicial inclusion. We refer to an inclusion $\varphi$ in $\Osc$ as a full inclusion of ordered simplicial complexes if $U\varphi$ is a full simplicial inclusion in $\sSetBC$.


     Now we would like to show that $U$ preserves pushouts along $Fj$. For if it does, then $p$ in the left diagram is taken to the pushout of $Z\leftarrow X \rightarrow UW$ on the right, which we know from the model structure on $\sSetBC$ (Theorem \ref{sSetBCmodel}) is a weak equivalence (as a pushout of $j \in J$):
     \[
     \begin{tikzcd} FX\arrow[r, "f"] \arrow[swap]{d}{Fj} & W \arrow[d, dashed, "p"] \\%
 FZ \arrow[r, dashed, "g"]& P
\end{tikzcd} \quad\quad \begin{tikzcd} UFX \cong X\arrow[r, "Uf"] \arrow[swap]{d}{j} & UW \arrow[d, dashed, "p'"] \\%
 UFZ \cong Z \arrow[r, dashed, "g'"]& P'
\end{tikzcd} \]
Given the pushout $P = FZ \xleftarrow{Fj} FX \xrightarrow{p} W$ in $\Osc$, we show $U(P)$ is the pushout of $UF(Z) \xleftarrow{U(j)} UF(X) \xrightarrow{U(f)} UW$ or equivalently the pushout of $Z \xleftarrow{j} X \xrightarrow{U(f)} UW$. Now, $UFj = j$ is a full simplicial map in $\sSetBC$. Recall that a $n$-loop is a sequence of $n$ directed 1-simplices starting and ending at the same vertex (see Definition \ref{loopdef}). It suffices to show that pushouts of $j$, $p': UW \rightarrow P'$, do not create any $n$-loops. Note that $FX$ is connected and thus only maps continuously to a single connected component of $W$, so we do not gain order relations in the pushout. Another way to see this is that to obtain the pushout $P'$, we follow the coequalizer construction as outlined in Lemma $\ref{osccolimits}$ with $f, Fj: FX \rightrightarrows W \sqcup FZ$. However, there are no added edges in the component $FZ$ since it has no extra order relations, so we take as in the construction $\colim(FX \rightrightarrows FZ \sqcup W')$, where $W'$ may have added edges. Since $FX \to FZ$ is a full simplicial inclusion, so is the map from $W'$ to the pushout, by Lemma \ref{lemma2}.  Therefore, vertices from $W'$ are not identified in the pushout, and the map $v$ on the discrete simplicial set of vertices in the diagram of Lemma $\ref{osccolimits}$ is the identity. Thus, the pushout $P$ is obtained by taking the vertices and simplices as those of $FP'$, and defining on the vertices the order relation given by $FZ \sqcup W$. Thus, if $P'$ does not have any $n$-loops, $F$ would not identify any simplices so that the simplices of $FP'$ would be the same as $P$. Then $UP \cong P'$ since the image under $U$ only depends on the simplices of $P$. Lemma \ref{oldpushout} gives this condition.

Thus, $P'$ has the same simplices as $P$, so that $U$ preserves pushouts along $Fj$ for $j \in J$. Now, maps in $cell(FJ)$ are the transfinite compositions of pushouts of the above form. Directed colimits are preserved by $U$, as established in Lemma $\ref{slemma2}$; this is since the resulting simplicial set has the same simplices, which are not affected by the loss of any extra order relations. A single pushout diagram in $\Osc$ with pushout $P$ gives a diagram whose pushout $P'$ in $\sSetBC$ is the simplicial set with the same simplices as $P$. Thus, $cell(J)$ are the weak equivalences in $\sSetBC$, and since each pushout of a map in $FJ$ is mapped to the corresponding pushout of $cell(J)$ in $\sSetBC$, any composition of pushout diagrams, by which we mean maps in $cell(FJ)$, is mapped to a composition of maps in $cell(J)$, which is a weak equivalence by the model structure on $\sSetBC$. Hence $U$ takes maps in $cell(FJ)$ to weak equivalences.
 \end{proof}
 
 Even though $(F, U)$ is not an equivalence of categories, we can use this adjunction to get an equivalence at the level of homotopy theories.
      \begin{theorem}
The model structure on $\sSetBC$ of Theorem \ref{sSetBCmodel} can be transferred to $\Osc$ via the adjunction  \begin{tikzcd}
            \sSetBC \arrow[r, shift left=.75ex,"F"{name=G}] & \Osc\arrow[l, shift left=.75ex, "U"{name=F}] 
            \arrow[phantom, from=F, to=G, "\scriptstyle\dashv" rotate=270]    
        \end{tikzcd} which therefore is a Quillen pair. Furthermore, it is a Quillen equivalence.
        \label{maintheorem2}
\end{theorem}
\begin{proof}
This is a again an application of Theorem \ref{thm1} for $(F, U)$ where $I$ is the set of generating cofibrations of $\sSetBC$ and $J$ is the set of generating trivial cofibrations of $\sSetBC$. $\Osc$ is closed under all small limits and colimits, by Lemma \ref{osccolimits}. The first condition is easily checked since a small simplicial set in $\sSetBC$ has a finite set of simplices by Lemma \ref{slemma}, and $U$ strictly decreases these simplices, so the small object argument is satisfied, by Lemma \ref{slemma2}. Lemma \ref{eq2} is the second condition. Hence $(F, U)$ is a Quillen pair.

 To see that it is a Quillen equivalence, note that due to the transferred model structure, $X \rightarrow Y$ is a weak equivalence in $\Osc$ precisely when $UX \rightarrow UY$ is a weak equivalence in $\sSetBC$. In order for $(F , U)$ to be a Quillen equivalence, we need that for all cofibrant $A$ in $\sSetBC$ and fibrant $B$ in $\Osc$, that $F(A) \rightarrow B$ is a weak equivalence in $\Osc$ precisely when $A \rightarrow U(B)$ is a weak equivalence in $\sSetBC$. Any cofibrant $A \in \sSetBC$ is weakly equivalent to $U(A')$ for some $A' \in \Osc$. This is because cofibrant objects of $\sSetBC$ are obtained via (retracts of) transfinite compositions of pushouts of the generating cofibrations of $\sSetBC$, starting with the pushout $\emptyset \sqcup_{\emptyset} \Delta[0]$. Lemma \ref{oldpushout} gives that pushouts of any generating cofibrations of $\sSetBC$ do not contain loops. Hence, $A \cong UF(A)$ for $A$ cofibrant, so we can take $A' \cong F(A)$. Suppose $F(A) \rightarrow B$ is a weak equivalence, then $UF(A) \cong A \rightarrow U(B)$ is a weak equivalence. Note that $FU(B)$ is weakly equivalent to $B$ since $UFU(B) \cong U(B)$ and the induced model structure. Then $UF(A) \rightarrow UFU(B)$ and hence $UF(A) \rightarrow U(B)$ is a weak equivalence, which implies that $F(A) \rightarrow B$ is a weak equivalence as required.
\end{proof}

\section{Properties of the Model Categories $\sSetBC$ and $\Osc$}

In this section, we will show certain properties of $\Osc$ and $\sSetBC$. In particular, we will give a partial classification of the cofibrant objects in $\sSetBC$ and $\Osc$, and show that their model structures are proper. All objects are cofibrant in $\sSet$, but the cofibrant objects in categories such as $\Cat$, $\Pos$, and $\sSetB$ are very difficult to understand, and we find this is the case for $\Osc$ and $\sSetBC$ also. Nevertheless, we identify some features of cofibrant objects in $\sSetBC$ and $\Osc$. We again refer to the generating trivial cofibrations of $\sSet$ as $I$, and trivial cofibrations as $J$.


 \begin{theorem}
    Cofibrant objects in $\sSetBC$ do not contain any loops in the sense of Definition \ref{loopdef}. In particular, any cofibrant object $X$ in $\sSetBC$ is of the form $X \cong UY$ for some ordered simplicial complex $Y \in \Osc$.
     
 \end{theorem}
 \begin{proof}
     The initial object in $\sSetBC$ is the empty simplicial set $\emptyset$, and the terminal object is $\Delta[0]$. 
     As stated in the proof of Theorem \ref{maintheorem2}, the cofibrant objects of $\sSetBC$ are the retracts of the transfinite compositions of pushouts of the generating cofibrations of $\sSetBC$, starting with the pushout $\emptyset \sqcup_\emptyset \Delta[0] = \Delta[0]$. We may iteratively build up other cofibrations starting from $\emptyset$ by composing with pushouts along generating trivial cofibrations in $\sSetBC$. Denote a transfinite composition of these pushouts by $\emptyset \to P$. By Lemma \ref{oldpushout}, no loops are formed along pushouts along the generating cofibrations in $\sSetBC$, so $P$ does not contain any loops. In addition, the cofibrant objects are the targets of the retracts of maps of the form $\emptyset \to P$, and therefore are simplicial subsets of $P$ (as the retracts must be levelwise inclusions of simplicial sets). Therefore, no cofibrant objects can contain any loops. 
     
     Since $F: \sSetBC \to \Osc$ was defined to collapse any loops in the simplicial set and then construct an ordered simplicial complex with the same simplices as the quotiented simplicial set (see Definition \ref{defquot}), for any cofibrant object $X\in\sSetBC$, we have that $FX$ is an ordered simplicial complex with the same simplices as $X$. Furthermore, $U$ gives a simplicial set with the same simplices as $FX$, and thus we know that for any cofibrant $X$ we have that $X \cong UFX$.
 \end{proof}

 To classify the cofibrant objects in $\Osc$, the following definition will be helpful.
 \begin{definition}
 An object in $\Osc$ is minimally ordered if its order relation is the transitive closure of the order relations imposed by its simplices.
 \end{definition}
 \begin{example}
 Consider the ordered simplicial complex given by the following diagram:
 \[
     \begin{tikzcd}
         x \leftarrow y \rightarrow x'
     \end{tikzcd}
      \]
      The ordered simplicial complex given by the order relation $y <x$ and $x < x'$ is minimally ordered. If we impose an order relation $x < x'$, it is no longer minimally ordered.
 \end{example}

 \begin{theorem}
     The cofibrant objects of $\Osc$ are those which can be written as the image $FX$ for a cofibrant object $X$ in $\sSetBC$. In particular, all cofibrant objects must be minimally ordered.
\end{theorem}
\begin{proof}
Cofibrant objects in $\Osc$ consist of the retracts of transfinite compositions of pushouts of the generating cofibrations in $\Osc$, which are the maps in the set $FD\Sd^2I$. First, we have that the initial object is the empty ordered simplicial complex, and the terminal object is the single vertex. The generating cofibrations, $FD\Sd^2I$, have the same simplices as the generating cofibrations in $\sSetBC$ (since the generating cofibrations of $\sSetBC$ do not contain any loops and thus no quotienting of vertices occurs for the involved simplicial sets when applying $F$).
    Now $D \Sd^2(\emptyset \hookrightarrow \Delta[0]) \cong (\emptyset \hookrightarrow *)$, so the cofibrant objects in $\Osc$ are obtained as transfinite compositions starting with $* \sqcup_\emptyset \emptyset = *$. We can iteratively build the pushout $P_n = FD\Sd^2\Delta[n] \sqcup_{FD\Sd^2 \partial\Delta[n]}P_{n-1}$, and suppose that $P_{n-1}$ also has no loops or no additional order relations and thus $P_{n-1} \cong FU(P_{n-1})$. 
    
    Now, let $W$ be any ordered simplicial complex. By the same argument as Lemma \ref{eq2}, in which the key point was that pushouts of the underlying simplicial sets involved did not have any loops, pushouts $FD\Sd^2\Delta[n] \sqcup_{F D \Sd^2 \partial\Delta[n]} W$ in $\Osc$ are mapped to $D\Sd^2\Delta[n] \sqcup_{D \Sd^2 \partial\Delta[n]} UW$ in $\sSetBC$. Therefore, we have that $UP_n \cong D\Sd^2\Delta[n] \sqcup_{D \Sd^2 \partial\Delta[n]}UP_{n-1}$. Applying $F$ to this pushout gives back $P_n$, since left adjoints commute with colimits and since $FU$ is the identity on each of the components of these particular pushouts. Denoting a transfinite composition as $P = FP'$ where $P'$ is a transfinite compositition of pushouts in $\sSetBC$, we obtain that in $\Osc$, cofibrant objects are simplicial subcomplexes of $P = FP'$, similar to the case in $\sSetBC$, and are obtained as $FX$ for some cofibrant $X$ in $\sSetBC$. 
    
    Note that by construction, any ordered simplicial complex in the image of $F$ is minimally ordered, meaning all cofibrant objects in $\Osc$ also are minimally ordered.
\end{proof}

Now we show that $\sSetBC$ and $\Osc$ are proper model categories. 

\begin{definition}[{\cite[Definition B.9]{HillHopkinsRavenel2016}}]
    A map $f \colon A \to X$ in a category $\mathcal{C}$ with small colimits is \textit{flat} if for every $A \to B$ and every weak equivalence $B \to B'$, the map $X \sqcup_A B \to X \sqcup_A B'$ is a weak equivalence.
\end{definition}

To prove left properness, by \cite[Example B10]{HillHopkinsRavenel2016}, it is sufficient to show that all generating cofibrations are flat, as flatness is preserved by pushouts and retracts. Thus we look at pushouts along generating cofibrations and weak equivalences.

\begin{lemma}
For any generating cofibration $\Sd^2\partial\Delta[n] \hookrightarrow  \Sd^2\Delta[n]$  and any weak equivalence $f\colon \Sd^2\partial\Delta[n] \to Y$, the map to the pushout $\Sd^2\Delta[n] \to \Sd^2\Delta[n]\sqcup_{\Sd^2\partial\Delta[n]}Y$ in $\sSetBC$ is a weak equivalence. 
   \label{leftproper}
\end{lemma}
\begin{proof}

    By the transferred model structure on $\sSetBC$ from $\sSetB$ (and so also from $\sSet$), a map $f: \Sd^2\partial\Delta[n] \to Y$ is a weak equivalence in $\sSetBC$ if and only if $\Ex^2f \simeq f$ is a weak equivalence in $\sSet$. Since $\sSet$ is proper, the pushout $\Sd^2\Delta[n] \to \Sd^2\Delta[n]\sqcup_{\Sd^2\partial\Delta[n]}Y$ is also a weak equivalence, and thus the geometric realization of $\Sd^2\Delta[n]\sqcup_{\Sd^2\partial\Delta[n]}Y$ is contractible. We would like to show that pushouts obtained in this way stay contractible once applying $LD$, which gives the pushout in $\sSetBC$. By \cite[Proposition 8.5]{fjellbonsSet}, $\sSetB$ is proper, so we only have to verify this after applying $L$.
    
    Similar to the proof of Lemma \ref{thmbc}, we have that the pushout in $\sSet$ can be viewed as attaching $Y$ to $\Sd^2\Delta[n]$ along the boundary inclusion. We can view the pushout as a contractible interior of simplices,  which are surrounded by possibly affected simplices, at most one edge away from the boundary, and lining the inside of the boundary. We need not worry about the simplices from $Y$, since in the geometric realization as they retract onto the boundary due to $Y$ being weakly equivalent to $\Sd^2\partial\Delta[n]$. Thus, we can restrict to considering only the quotiented boundary of $\Sd^2\Delta[n]$ and its interior, which together we will denote as $P$. 
    
    Note $P$ is contractible as the pushout of a weak equivalence along a cofibration to a contractible simplicial set in $\sSet$. We then take a deformation retraction in the geometric realization of $P$ to its barycenter, which lies within the the contractible interior; this is a family of paths starting from the boundary. Now, $L$ may collapse pairs of simplices within the affected simplex area; but importantly, any pair of parallel simplices are faces of a pair of parallel $n$-simplices, due to the subdivision, so identifying the appropriate paths reaching such parallel faces would still result in a well-defined deformation retraction. To elaborate, after applying $L$, the parallel faces would be identified as they are collapsed to the same face, but so would the $n$-simplices which they are faces of, hence there is a unique path for the resultant deformation retraction; this is similar the process in Lemma \ref{defret2}. Thus $LP$ is contractible, so the induced map $\Sd^2\Delta[n]\sqcup_{\Sd^2\partial\Delta[n]}Y \to LD(\Sd^2\Delta[n]\sqcup_{\Sd^2\partial\Delta[n]}Y)$ is a weak equivalence, inducing a map of contractible spaces in the geometric realization.  
\end{proof}

\begin{theorem}
   The model structure on $\sSetBC$ from Theorem \ref{sSetBCmodel} is proper.
    \label{BCproper}
\end{theorem}
\begin{proof}
    Since $\sSetBC$ is obtained via right transfer (see Theorem \ref{sSetBCmodel}) from the proper model category $\sSetB$ \cite[Theorem 1.2]{fjellbonsSet}, it is right proper. By Lemma \ref{leftproper}, and \cite[Proposition B.11]{HillHopkinsRavenel2016}, all cofibrations are flat since the generating ones are, so $\sSetBC$ is left proper, and hence proper.
\end{proof}
\begin{theorem}
    The model structure on $\Osc$ from Theorem \ref{maintheorem2} is proper.
\end{theorem}
\begin{proof}
    $\Osc$ is right proper since it is obtained via right transfer from $\sSetBC$, which is proper by Theorem \ref{BCproper}. To show that $\Osc$ is left proper, we again show that generating cofibrations are flat. This is the case due to Lemma \ref{oldpushout} and the proof of Theorem \ref{maintheorem2}, since loops cannot be formed along generating cofibrations in the underlying pushout diagrams in $\sSetBC$. Thus, by the construction of coequalizers (and hence pushouts) in Lemma \ref{osccolimits}, these particular pushouts in $\Osc$ are precisely the same as applying $F$ to the pushouts of the underlying diagrams in $\sSetBC$. Then, applying $U$ to the pushout in $\Osc$ gives back the pushout of the underlying diagram in $\sSetBC$. 
    
    A pushout along a generating cofibration must then preserve weak equivalences in $\Osc$. To see this, denote the map to the pushout as $p\colon F(\Sd^2\Delta[n]) \to F(\Sd^2\Delta[n]) \sqcup_{F(\Sd^2\partial\Delta[n])} Y$ in $\Osc$, where $F(\Sd^2\partial\Delta[n])\to Y$ is a weak equivalence. Let $P := F(\Sd^2\Delta[n]) \sqcup_{F(\Sd^2\partial\Delta[n])} Y$. By the previous discussion, we have that $UP \cong \Sd^2\partial\Delta[n] \sqcup_{\Sd^2\partial\Delta[n]} UY$, the pushout of the underlying diagram in $\sSetBC$. Applying $U$ to the pushout diagram, we get $Up: \Sd^2\Delta[n] \to \Sd^2\partial\Delta[n] \sqcup_{\partial\Sd^2\Delta[n]} UY$ in $\sSetBC$, where $\Sd^2\partial\Delta[n] \to UY$ is still a weak equivalence.  Since $\sSetBC$ is left proper (Theorem \ref{BCproper}), the map $Up:\Sd^2\partial\Delta[n] \to UP$ is a weak equivalence; but then so is $p \colon F(\Sd^2\partial\Delta[n]) \to F(UP)$ in $\Osc$ (note that in this instance we have $F(UP) \cong P$ and $UF(\Sd^2\partial\Delta[n]) \cong \Sd^2\partial\Delta[n] $). Then since the generating cofibrations are flat, all cofibrations are flat, and thus $\Osc$ is left proper.
\end{proof}

  \bibliographystyle{alpha}

 \end{document}